\documentclass[a4paper,oneside,11pt]{article}

\usepackage{amsmath}
\usepackage{amssymb}
\usepackage[utf8]{inputenc} %try and remove if possible
\usepackage{accents}
\usepackage{mathtools}
\usepackage{amsthm}
\usepackage[english]{babel}
\usepackage{letltxmacro}
\usepackage{xparse}
\usepackage[dvipsnames,svgnames,x11names,hyperref]{xcolor}
\usepackage{graphicx}

\usepackage[backgroundcolor = white, linecolor = Cerulean, textcolor = Cerulean, bordercolor = Cerulean]{todonotes}

\usepackage{fullpage}
\usepackage{fancyhdr}
\usepackage{ebgaramond}
\usepackage{euscript}
\usepackage{mathrsfs}
\usepackage{accents}
\usepackage[bbgreekl]{mathbbol}
\usepackage[euler]{textgreek}
\usepackage[scr=boondoxo]{mathalfa}
\DeclareSymbolFontAlphabet{\mathbbm}{bbold}
\DeclareSymbolFontAlphabet{\mathbb}{AMSb}%
\usepackage{mleftright}
\usepackage{dirtytalk}
\usepackage{shuffle}
\usepackage{scalerel}
\usepackage{tensor}
\usepackage{soul}
\setul{0pt}{0.4pt}

\makeatletter
\newcommand*\bigcdot{\mathpalette\bigcdot@{.5}}
\newcommand*\bigcdot@[2]{\mathbin{\vcenter{\hbox{\scalebox{#2}{$\m@th#1\bullet$}}}}}
\makeatother

\mathchardef\mhyphen="2D

\DeclareFontFamily{U}{mathx}{\hyphenchar\font45}
\DeclareFontShape{U}{mathx}{m}{n}{
	<5> <6> <7> <8> <9> <10>
	<10.95> <12> <14.4> <17.28> <20.74> <24.88>
	mathx10
}{}
\DeclareSymbolFont{mathx}{U}{mathx}{m}{n}
\DeclareFontSubstitution{U}{mathx}{m}{n}
\DeclareMathAccent{\widecheck}{0}{mathx}{"71}
\DeclareMathAccent{\wideparen}{0}{mathx}{"75}

\usepackage[ruled]{algorithm2e}

\usepackage{listings}
\usepackage{tikz}
\usetikzlibrary{decorations.pathreplacing}
\usetikzlibrary{decorations.markings}
\usepackage{pgfplots}
\usepackage{wrapfig}
\usepackage{caption}
\usepackage[percent]{overpic}
\usepackage{tikz-cd}

\usepackage{aliascnt}
\usepackage[pdftex, ocgcolorlinks, colorlinks = true,
linkcolor = imperialBrick,
citecolor = imperialGreen,
urlcolor = imperialProcess
]{hyperref}

\usepackage{nameref}

\addto\extrasenglish{
	
} % This makes Chapter appear with a capital when referenced

\addto\extrasenglish{
	
} % This makes Section appear with a capital when referenced

\addto\extrasenglish{
	
} % This makes Section appear with a capital when referenced

\newcommand{\mynewtheorem}[2]{
	\newaliascnt{#1}{dummy}
	\newtheorem{#1}[#1]{#2}
	\aliascntresetthe{#1}
	\expandafter\def\csname #1autorefname\endcsname{#2}
}

\newtheoremstyle{note}
{\topsep}   % ABOVESPACE
{\topsep}   % BELOWSPACE
{}  % BODYFONT
{0pt}       % INDENT (empty value is the same as 0pt)
{\itshape\bfseries} % HEADFONT
{.}         % HEADPUNCT
{5pt plus 1pt minus 1pt} % HEADSPACE
{{\color{imperialTangerine}\thmname{#1}\emph{\thmnumber{ #2}}}\hspace{0.1em}\textnormal{\thmnote{ (#3)}}}     % CUSTOM-HEAD-SPEC

\theoremstyle{plain}
\mynewtheorem{thm}{Theorem}
\mynewtheorem{thmdef}{Theorem/Definition}
\mynewtheorem{prop}{Proposition}
\mynewtheorem{cor}{Corollary}
\mynewtheorem{lem}{Lemma}
\mynewtheorem{conj}{Conjecture}
\theoremstyle{definition}
\mynewtheorem{defn}{Definition}
\mynewtheorem{expl}{Example}
\mynewtheorem{prob}{Problem}
\mynewtheorem{ques}{Question}
\mynewtheorem{cond}{Condition}
\mynewtheorem{conv}{Convention}
\mynewtheorem{defthm}{Definition/Theorem}
\theoremstyle{remark}
\mynewtheorem{rem}{Remark}

\theoremstyle{note}

\definecolor{imperialBlue}{RGB}{0, 62, 116}
\definecolor{imperialBrick}{RGB}{165,25,0}
\definecolor{imperialProcess}{RGB}{0,133,202}
\definecolor{imperialGreen}{RGB}{2,137,59}
\definecolor{imperialRed}{RGB}{221,37,1}
\definecolor{imperialOrange}{RGB}{210,64,0}
\definecolor{imperialBlue2}{RGB}{0,110,175}
\definecolor{imperialTangerine}{RGB}{236,115,0}
\definecolor{imperialPurple}{RGB}{101,48,152}
\definecolor{imperialLime}{RGB}{196,214,0}
\definecolor{imperialKermit}{RGB}{102,164,10}

\definecolor{kermit}{RGB}{102,164,10}
\definecolor{teal}{RGB}{0,142,170}
\definecolor{tangerine}{RGB}{236,115,0}
\definecolor{raspberry}{RGB}{145,0,72}
\definecolor{lime}{RGB}{196,214,0}

\colorlet{chaptergrey}{imperialBlue}

\usepackage{cutwin,kantlipsum,mwe}
\usepackage{regexpatch}
\opencutright

\makeatletter
\xpretocmd{\cutout}{\leavevmode\hrule \@height\z@ \@width\linewidth\relax}
\makeatother

\usepackage{iftex}
\iftutex
\usepackage{unicode-math}
\setmathfontface\altgrfont{GFS Artemisia Italic}[Scale=MatchLowercase]
\newcommand{\text{\textdelta}}{\altgrfont{δ}}
\else
\usepackage[T1]{fontenc}
\DeclareSymbolFont{altgr}{OML}{antt}{m}{it}
\DeclareMathSymbol{\sko }{\mathord}{altgr}{"0E}
\fi

\let\oldKronecker\sko
\renewcommand{\sko}{\oldKronecker \hspace{-0.1em}}

\def\skoo{\sko\hspace{0.1em}}
\def\kron{\text{\textdelta}}

\makeatletter
\def\upintkern@{\mkern-7mu\mathchoice{\mkern-3.5mu}{}{}{}}
\def\upintdots@{\mathchoice{\mkern-4mu\@cdots\mkern-4mu}%
	{{\cdotp}\mkern1.5mu{\cdotp}\mkern1.5mu{\cdotp}}%
	{{\cdotp}\mkern1mu{\cdotp}\mkern1mu{\cdotp}}%
	{{\cdotp}\mkern1mu{\cdotp}\mkern1mu{\cdotp}}}

\newcommand{\UpMultiIntegral}[1]{%
	\edef\ints@c{\noexpand\upintop
		\ifnum#1=\z@\noexpand\upintdots@\else\noexpand\upintkern@\fi
		\ifnum#1>\tw@\noexpand\upintop\noexpand\upintkern@\fi
		\ifnum#1>\thr@@\noexpand\upintop\noexpand\upintkern@\fi
		\noexpand\upintop
		\noexpand\ilimits@
	}%
	\futurelet\@let@token\ints@a
}
\makeatother

\DeclareFontFamily{OMX}{mdbch}{}
\DeclareFontShape{OMX}{mdbch}{m}{n}{ <->s * [0.8]  mdbchr7v }{}
\DeclareFontShape{OMX}{mdbch}{b}{n}{ <->s * [0.8]  mdbchb7v }{}
\DeclareFontShape{OMX}{mdbch}{bx}{n}{<->ssub * mdbch/b/n}{}

\DeclareSymbolFont{uplargesymbols}{OMX}{mdbch}{m}{n}
\SetSymbolFont{uplargesymbols}{bold}{OMX}{mdbch}{b}{n}
\DeclareMathSymbol{\upintop}{\mathop}{uplargesymbols}{82}
\DeclareMathSymbol{\upointop}{\mathop}{uplargesymbols}{"48}

\DeclareFontEncoding{MDB}{}{}
\DeclareFontFamily{MDB}{mdbch}{}
\DeclareFontShape{MDB}{mdbch}{m}{n}{ <->s * [0.8]  mdbchrmb }{}
\DeclareFontShape{MDB}{mdbch}{b}{n}{ <->s * [0.8]  mdbchbmb }{}
\DeclareFontShape{MDB}{mdbch}{bx}{n}{<->ssub * mdbch/b/n}{}
\DeclareFontSubstitution{MDB}{cmr}{m}{n}
\DeclareSymbolFont{mathdesignB}{MDB}{mdbch}{m}{n}%
\SetSymbolFont{mathdesignB}{bold}{MDB}{mdbch}{b}{n}%
\DeclareMathSymbol{\upintclockwise}{\mathop}{mathdesignB}{128}
\DeclareMathSymbol{\upointclockwise}{\mathop}{mathdesignB}{130}
\DeclareMathSymbol{\upointctrclockwise}{\mathop}{mathdesignB}{132}
\DeclareMathSymbol{\upoiint}{\mathop}{mathdesignB}{134}
\DeclareMathSymbol{\upoiiint}{\mathop}{mathdesignB}{136}

\makeatletter
\newcommand{\upint}{\DOTSI\upintop\ilimits@}
\newcommand{\upoint}{\DOTSI\upointop\ilimits@}
\makeatother

\renewcommand{\int}{\upint}

\usepackage[scr=boondoxo]{mathalfa}

\newcommand{\bbE}{\mathbb E}
\newcommand{\bbR}{\mathbb R}

\usepackage[euler]{textgreek}

\newcommand{\dif}{\mathrm{d}}

\def\EuD{\EuScript{D}}
\def\EuS{\EuScript{S}}

\def\bfX{\boldsymbol{X}}

\newcommand{\p}{{\lfloor p \rfloor}}

\usepackage{authblk}

\title{\textsc{On the Wiener Chaos Expansion of the\\ Signature of a Gaussian Process}}

\author[1]{Thomas Cass}

\author[2]{Emilio Ferrucci\thanks{Corresponding author: \href{mailto:Emilio.RossiFerrucci@maths.ox.ac.uk}{\texttt{Emilio.RossiFerrucci@maths.ox.ac.uk}}}}

\affil[1]{\small Department of Mathematics, Imperial College London}

\affil[2]{\small Mathematical Institute, University of Oxford}

\date{\today}

\makeatletter
\newcommand*{\starsection}[1]{%
	\section*{#1}%
	\NR@gettitle{#1}%
}
\makeatother

\usepackage{stackengine}

\usetikzlibrary{decorations.markings}
\usepackage[toc,page]{appendix}

\begin{document}
	
	\maketitle	
	\begin{abstract}
		We compute the Wiener chaos decomposition of the signature for a class of Gaussian processes, which contains fractional Brownian motion (fBm) with Hurst parameter $H \in (1/4,1)$. At level $0$, our result yields an expression for the expected signature of such processes, which determines their law \cite{CL13}. In particular, this formula simultaneously extends both the one for $1/2 < H$-fBm \cite{Bau07} and the one for Brownian motion ($H = 1/2$) \cite{Faw04}, to the general case $H > 1/4$, thereby resolving an established open problem. Other processes studied include continuous and centred Gaussian semimartingales.
	\end{abstract}
	
	\subsection*{Acknowledgements} The research of both authors is currently supported by EPSRC Programme Grant EP/S026347/1. The research of Emilio Ferrucci, while at Imperial College London, was supported by the Centre for Doctoral Training in Financial Computing \& Analytics EP/L015129/1 and later by the Strategic Project Grant EP/W522673/1.
	
	\section*{Introduction}
	
	The signature of a path $X \colon [0,T] \to \bbR^d$,
	\begin{equation}\label{eq:sig}
		\EuS(X)_{0T} \coloneqq \sum_{n = 0}^\infty\int_{0 < u_1 < \ldots < u_n < T} \dif X_{u_1} \otimes \cdots \otimes \dif X_{u_n} \quad \in T(\!(\bbR^d)\!)\text{ ,}
	\end{equation}
	is a series of tensors which, up to \say{retracings}, determines the image of $X$ \cite{HL10, BGLY16}. The probabilistic counterpart to this result states that, in many cases of interest, the law of a stochastic process is determined by its expected signature \cite{CL13}, which is therefore seen to play a role for processes analogous to that of moments for random variables. 
	
	The best-known example of an explicit formula for the expected signature of a stochastic process occurs in the case of Brownian motion: calling $\{e_1,\ldots,e_d\}$ the canonical basis of $\bbR^d$, we have
	\begin{equation}\label{eq:sumSquares}
		\mathbb E \EuS(X)_{st} = \exp\bigg( \frac {t-s}2 \sum_{\gamma = 1}^d e_\gamma^{\otimes 2}\bigg) = \sum_{n = 0}^\infty \frac{(t-s)^n}{2^n n!} \sum_{\gamma_1,\ldots,\gamma_n = 1}^d e_{\gamma_1}^{\otimes 2} \otimes \cdots \otimes e_{\gamma_n}^{\otimes 2}\text{ .}
	\end{equation}
	This identity was first shown by \cite{Faw04,LV07cub}, and later proved in a variety of different ways \cite{Bau04,FS17}. The expected signature of Brownian motion has also been studied in the case in which the process is stopped upon hitting the boundary of a domain \cite{LN15,BDMN21,LN22}.
	
	In \cite{Bau07} the authors derive an integral expression for the expected signature of fractional Brownian motion (fBm) with Hurst parameter $H \in (1/2,1)$. This result was extended in \cite{BPQ13,horatioPhd} to a more general class of Gaussian Volterra processes with sample paths that are more regular than Brownian motion, with the formula for the expected signature written in terms of the Volterra kernel. The method used involves a piecewise-linear interpolation of the paths of the process $X$, which reduces the calculation to that of a sum of mixed Gaussian moments, to which Wick's theorem applies, followed by a convergence argument. The expression in \cite{Bau07} does not, however, yield the correct prediction for the case of Brownian motion $H = 1/2$. When $H < 1/2$ it involves integrals that do not converge at all, and new ideas are needed to obtain a formula. On a technical level, the reason for these differences can be seen by considering the expression for the expected signature of a scalar $1/2 < H$-fBm $X$ at level $2$: calling $R(s,t) \coloneqq \mathbb E[X_sX_t]$ the covariance function of $X$, the formula states that
	\begin{equation}\label{eq:sigfBm2}
		\mathbb E \EuS(X)^{(2)}_{st} = \int_{s<u<v<t} R(\dif u, \dif v) = H(2H-1) \int_{s<u<v<t} (v-u)^{2H-2} \dif u \dif v\text{ .}
	\end{equation}
	Integrating either of the two variables generates an evaluation $(v-u)^{2H-1}|_{u = v}$, which is only finite when $H > 1/2$ and indeterminate when $H = 1/2$. In fact, approximating $X$ with a sequence of piecewise linear processes $(X^\ell)_{\ell \in \mathbb N}$ one obtains a sequence of integrals (actually finite sums) $\int_{s < u < v < t} \bbE[\dot X_u^\ell \dot X^\ell_v] \dif u \dif v$ which converges to the above double integral when $H > 1/2$, to $(t-s)/2$ when $H = 1/2$ (as predicted by \eqref{eq:sumSquares}), and continues to converge to $(t-s)^{2H}/2$ for $1/4 < H \leq 1/2$. When $H \leq 1/4$ the iterated integrals (in particular the Lévy area) of smooth approximations of $X$ do not converge in mean square, and other techniques (e.g.\ \cite{NuaTin11}) must be relied upon to define a rough path, and hence a signature. These rough paths present a number of differences with the canonical one defined for $H > 1/4$, and are therefore not considered in this paper.
	\begin{figure}[h]
		\minipage{0.5\textwidth}
		\includegraphics[width=\linewidth]{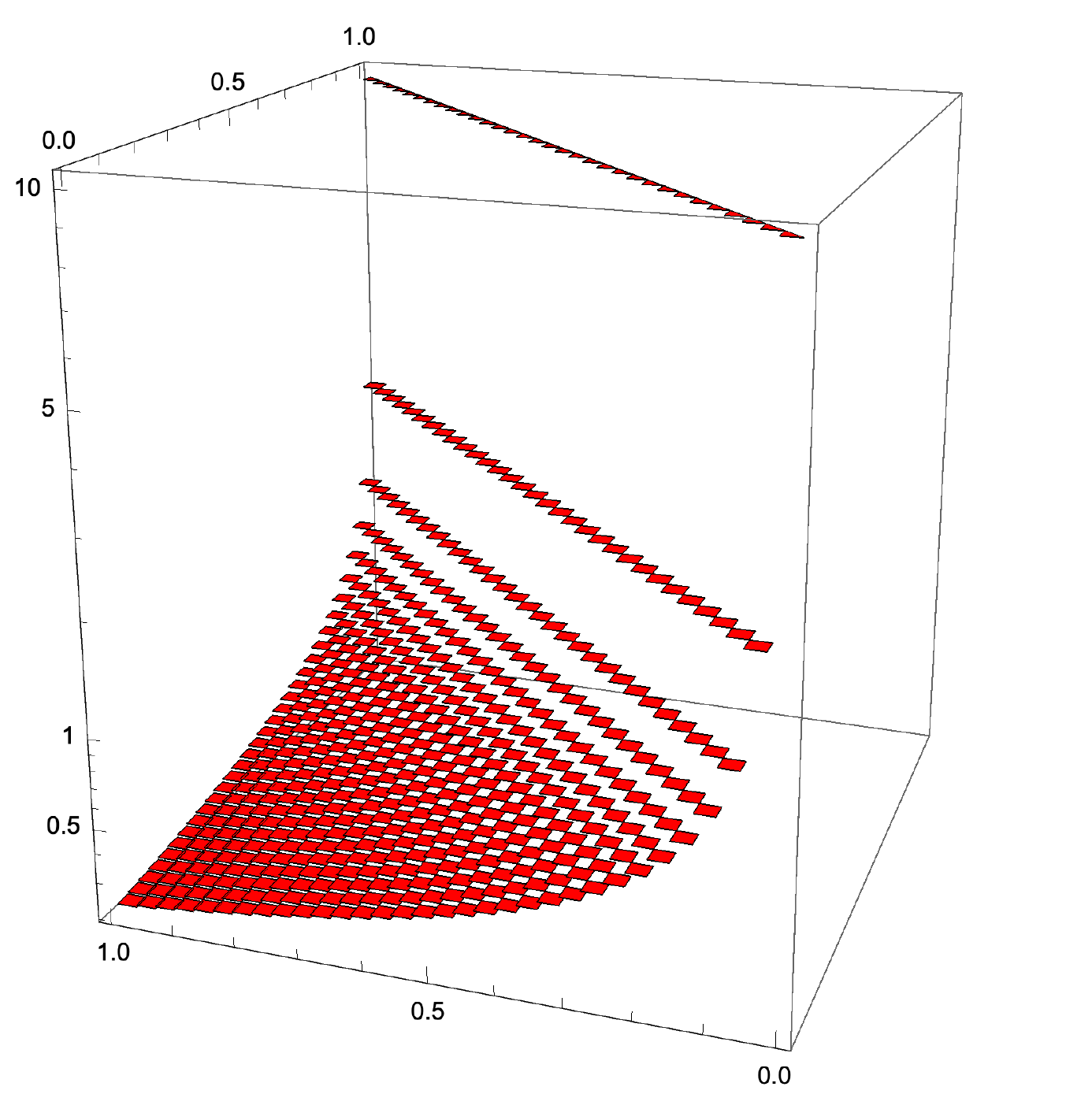}
		\vspace{-20pt}
		\endminipage\hfill
		\minipage{0.5\textwidth}
		\includegraphics[width=\linewidth]{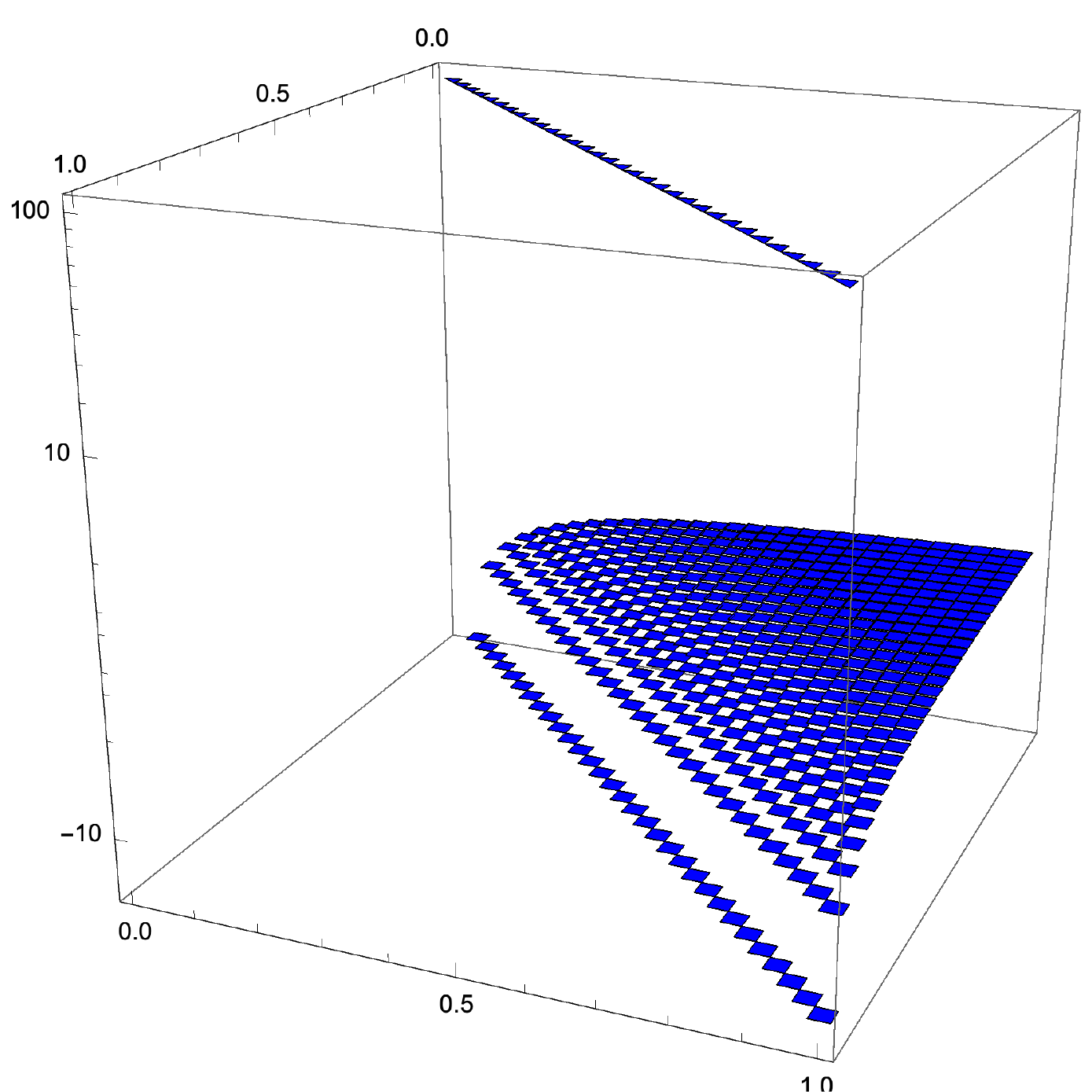}
		\vspace{-20pt}
		\endminipage
		\caption{Here we compare the two behaviours, corresponding to $H >1/2$ and $H<1/2$, of $\int_{0 < u < v < 1} \bbE[\dot X_u^\ell \dot X^\ell_v] \dif u \dif v$ with $X^\ell$ the sequence piecewise linear interpolations of $X$ on a partition. On the left we have chosen $H = 2/3$, and the sequence of integrals converges to a finite improper integral, whereas on the right $H = 1/3$ and the on- and off-diagonal contributions diverge to opposite infinities. (The plots are oriented in different ways and the $z$-axis is rescaled, both for improved visibility.) This graphic has been created using Wolfram Mathematica.}\label{fig:plusMinusCorr}
	\end{figure}
	
	What is needed to obtain a formula for the expected signature that also works in the case of negatively-correlated increments $1/4 < H < 1/2$ is a way of expressing the indeterminacy \say{$\infty - \infty$} explained in \autoref{fig:plusMinusCorr}. The trick for doing this is simple to describe: integrate out the first variable in \eqref{eq:sigfBm2} and, calling $R(t) \coloneqq R(t,t)$ the variance function of $X$, note that for $H > 1/2$ we have
	\begin{equation}\label{eq:keySub}
		\int_{s<u<v<t} \partial_{12} R(u, v) \dif u \dif v = \int_s^t \big[ \partial_2 R(v, v) - \partial_2 R(s,v) \big] \dif v = \int_s^t \big[ \tfrac 12 R'(v) - \partial_2 R(s,v) \big] \dif v \text{ .}
	\end{equation}
	We have replaced $\partial_2 R(v, v)$ with $\frac 12 R'(v)$, which can be done by symmetry of $R$:
	\begin{equation}\label{eq:RSub}
		R'(v) = \frac{\dif}{\dif v} R(v,v) = \partial_1 R(v,v) + \partial_2 R(v,v) = 2\partial_2 R(v,v) \text{ .}
	\end{equation}
	This is relevant to the case of $(1/4,1/2) \ni H$-fBm since, while $\partial_2 R(v,v)$ or $\partial_1 R(v,v)$ is the infinite evaluation discussed earlier, the last integral in \eqref{eq:keySub} is perfectly well defined. These integrands can be chained together on simplices, e.g.\ $\int_{s < u < v < t} [\tfrac 12 R'(u) - \partial_2 R(s,u)][\tfrac 12 R'(v) - \partial_2 R(u,v)] \dif u \dif v$, and combined with the other types of integrand $\partial_{12} R(w, z)$, to yield a formula that is very similar to that of \cite{Bau07}, but continues to be convergent for $1/4 < H < 1/2$ and agrees with \eqref{eq:sumSquares} for $H = 1/2$.
	
	Showing that the formula obtained by such substitution actually coincides with the expected signature for $X$ in a broad class of Gaussian processes --- essentially those Gaussian rough paths introduced in \cite{CQ02,LQ02,FV10b} with the imposition of a few additional smoothness and regularity requirements on the (co)variance function --- is the main focus of this paper. In fact, our main result will prove a formula for the full Wiener chaos expansion of $\EuS(X)$, the $0^\text{th}$ level of which is the expectation. As far as we know, the expression for the positive chaos projections of the signature is not to be found in the literature even in the classical case of Brownian motion. While the expression of the positive levels of Wiener chaos is very similar in spirit to that of the $0^\text{th}$, it requires us to use some Malliavin calculus in the setting of $1$-parameter Gaussian processes, and results in technical complications in the proof of convergence. The main additional ingredients needed are Stroock's formula for the $m$-th Wiener chaos projection and a novel definition of multiple Wiener integral of a function. For the latter, it should be noted that while multivariate, deterministic integrands for Gaussian noise naturally live in a certain Hilbert space (which for fBm can be identified with a Sobolev space), we are interested in integrating functions of multiple times, i.e.\ $\int_{[0,T]^m} f(t_1,\ldots,t_m) \dif X^{\gamma_1}_{t_1} \cdots \dif X^{\gamma_m}_{t_m}$ in a Skorokhod-type sense: this is achieved by approximating $f$ with elementary integrands, and showing independence of the approximation. Computing the Wiener chaos projections of the signature of a Gaussian process $X$ has the benefit of expressing $\EuS(X)$ as a sum of terms that are orthogonal in $L^2$, something that has the potential to be used for various types of numerical calculations, e.g.\ estimates of Euler expansions for Gaussian rough differential equations. It should be mentioned that, while (in the cases considered) the expected signature already determines the law of $X$ and therefore that of the Wiener chaos projections of $\EuS(X)$, it does not appear obvious how one may obtain the latter from the former directly. While fBm is the main example of a process for which our calculation is novel, we briefly also consider centred, continuous Gaussian semimartingales, such as the Brownian bridge returning to the origin and centred Ornstein-Uhlenbeck processes with deterministic initial condition.
	
\begin{figure}[h]
		\minipage{0.45\textwidth}
		\includegraphics[width=\linewidth]{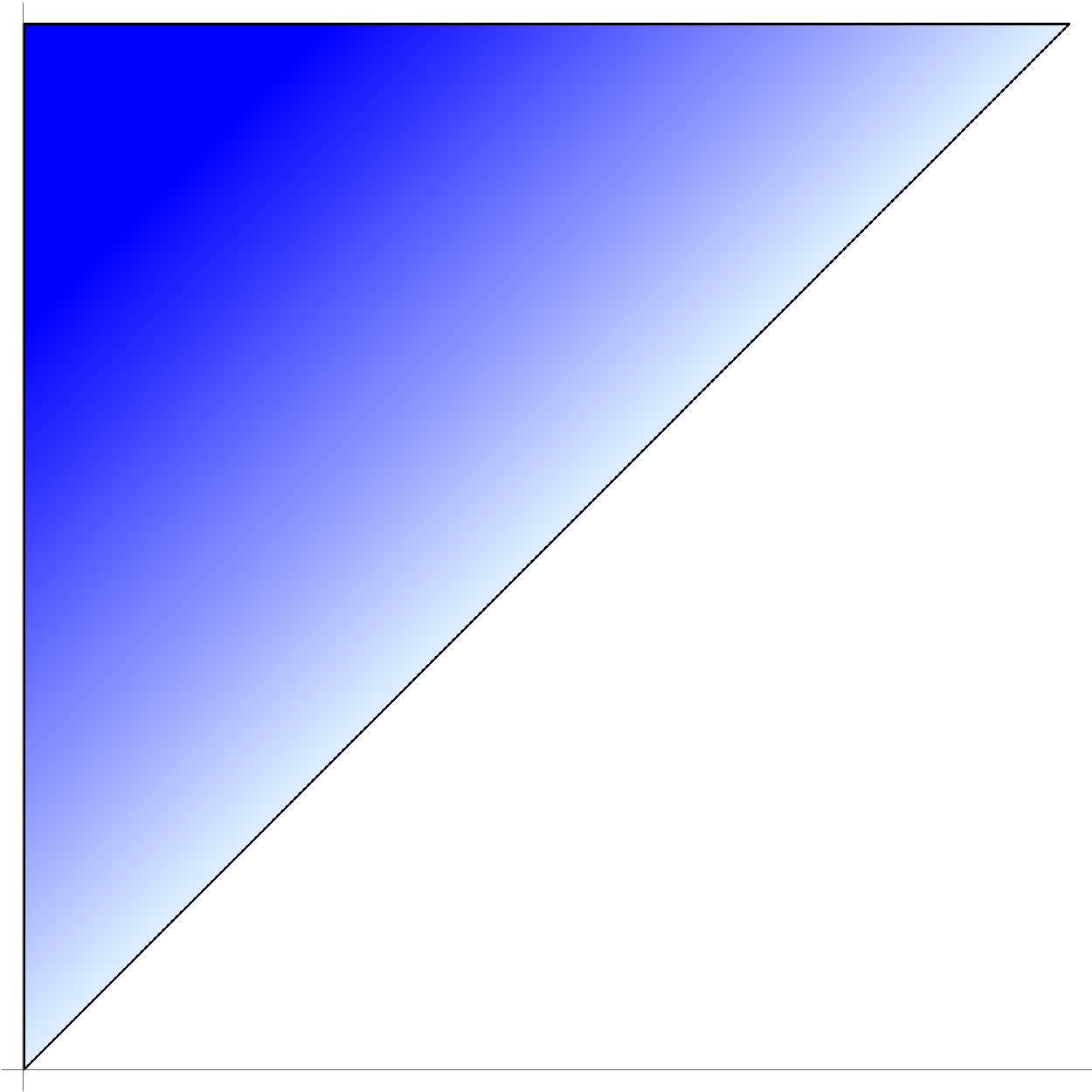}
		\vspace{-20pt}
		\endminipage\hfill
		\minipage{0.45\textwidth}
		\includegraphics[width=\linewidth]{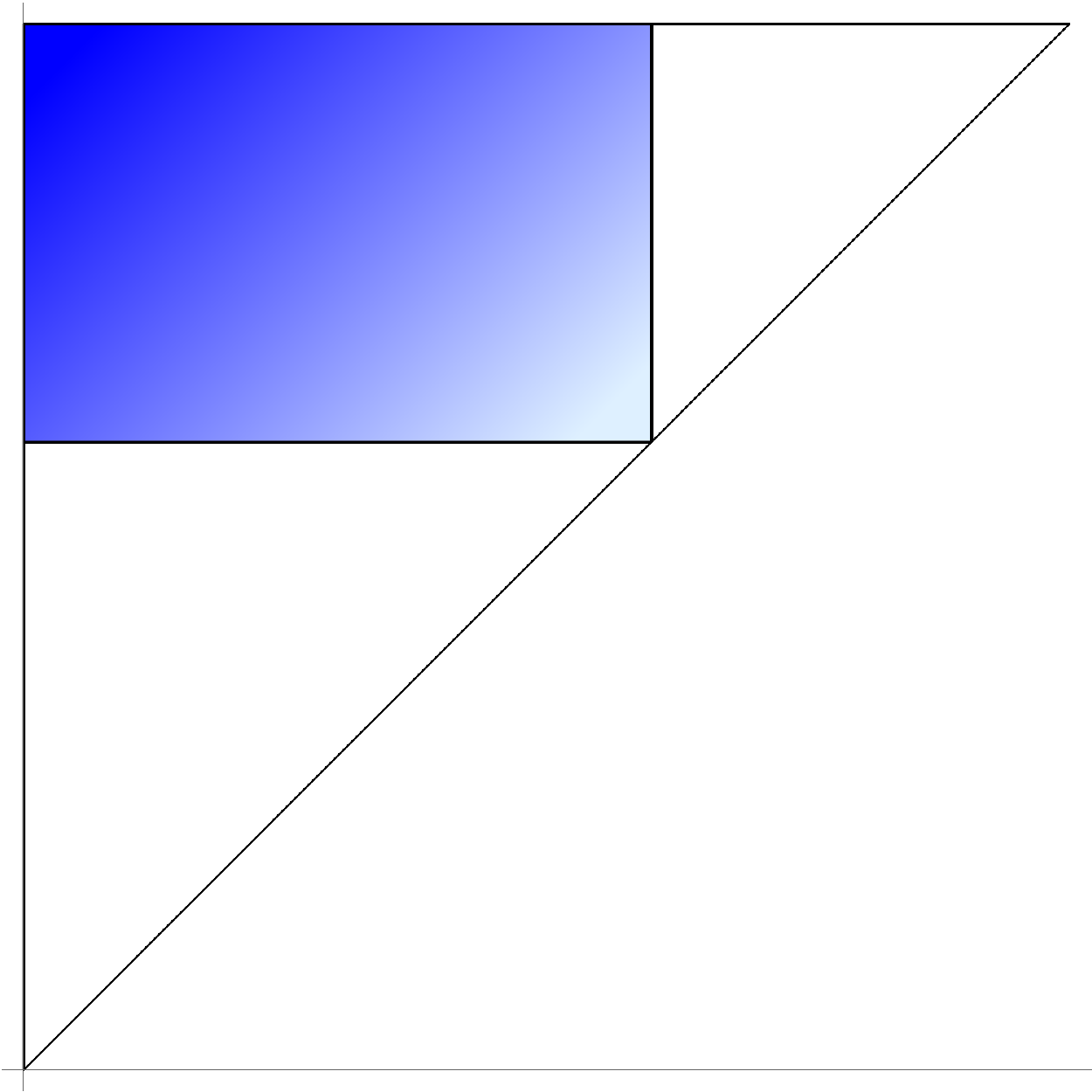}
		\vspace{-20pt}
		\endminipage
		\caption{A graphic representing the contour plot of $(t-s)^{2H-2}$ on $\{0 < s < t < 1\}$ (on the left) and $\{0 < s < u < t < 1\}$ with $u \in [0,1]$ fixed (on the right): the integral of the former is improper on the whole diagonal, while that of the later only at a point: when $0 < H < 1/2$, only the latter converges. This graphic has been created using Wolfram Mathematica.}\label{fig:gradient}
	\end{figure}

	As in the main reference article \cite{Bau07}, the technique that underlies our proof is piecewise-linear approximation of $X$. The arguments needed to prove the result are however much more involved, for three essential reasons. First is the fact that we must perform and justify the substitution \eqref{eq:keySub}, which requires novel arguments for convergence; even proving finiteness of the integrals in the main formula requires more sophisticated bounds in the $1/4 < H < 1/2$ case than it does in the $H > 1/2$ case (see \autoref{fig:gradient} for the simplest example of an observation that must be made when $H < 1/2$). Second is that Malliavin derivatives are involved for positive levels of the Wiener chaos and third is that our arguments must accommodate a wider class of Gaussian processes.
	While the substitution \eqref{eq:keySub} may seem very natural, it does not emerge obviously from the proof that we have given here, and must instead be guessed in advance. Indeed, it is worth mentioning that the way in which we first derived the statement of the main result involved an entirely different approach, which made use of the Skorokhod-rough integral conversion formula \cite{CL19,CL20}, applied recursively to the RDE for the signature. The outline of this proof can be found in the second named author's PhD thesis \cite[Ch.\ 5]{Fer22}. While this approach has the drawback of generating further technical problems, reason for which it is not the one presented here, it has the advantage of leading up constructively to the main formula.
	
	This paper is organised as follows: in \autoref{sec:back} we briefly introduce the class of Gaussian processes considered and the Malliavin calculus framework for them; we then use this language to identify functions as multiple Wiener integrands. In \autoref{sec:main} we state the main result \autoref{thm:main} and discuss a few consequences and examples that follow; in \autoref{sec:proof} we prove the main result; in \nameref{sec:concl} we outline some aspects that could be tackled in further research. Finally, it should be mentioned that in \cite{Bau07}, in addition to the expected signature of $1/2 < H$-fBm, the authors also compute the expected signature at levels $2$ and $4$ for $1/4 < H$-fBm in a manner that does not obviously generalise to different processes or higher levels; while not necessary in our proofs, it is sensible to verify that our main result agrees with this calculation: this check is performed in \autoref{appendix}.

	\section{Background on Malliavin calculus for Gaussian processes}\label{sec:back}
	
	In this section we introduce the class of Gaussian processes to which this paper applies, establish some notation, and give a brief overview of the tools of Malliavin calculus that are necessary in the proof of the main result. We follow \cite{Nua06,NouPec12} for the general Malliavin calculus framework, \cite{HuCa78} for its aspects that pertain to Gaussian processes indexed by a time parameter, and \cite{CF10,CL19,CL20} for aspects regarding the rough path lifts of such processes.
	
	Throughout this paper we will be working with a Gaussian process with i.i.d.\ components $X \colon \Omega \times [0,T] \to \bbR^d$ where $\Omega = C([0,T],\bbR^d)$, $X_t(\omega) \coloneqq \omega(t)$, $\mathcal F_t \coloneqq \sigma(X_{s} : 0 \leq s \leq t)$. We assume $X$ to be centred, i.e.\ $\mathbb E X \equiv 0$, and for it to have deterministic initial condition $X_0 = 0$. We will write $X_{st} \coloneqq X_t - X_s$ for the increments of $X$. By Gaussianity, the probability measure $\mathbb P$ on $\Omega$ is characterised by the covariance function of $X$
	\begin{equation}
		R \colon [0,T]^2 \to \bbR^d \otimes \bbR^d, \quad R(s,t) \coloneqq \bbE[X_s \otimes X_t]\text{ .}
	\end{equation}
	We will denote $R(\,\cdot\,)$ the variance function of $X$, i.e.\ $R(t) \coloneqq R(t,t)$. The independence hypothesis implies that $R$ is a diagonal matrix $R^{\alpha\beta} = \kron^{\alpha\beta} R^{\alpha\alpha}$, and the fact that they are identically distributed, $R^{\alpha\beta} = \kron^{\alpha\beta}R^{11}$ will be determined by a single scalar function, which by abuse of notation we will also call $R$. Although our results can be conjectured to continue to hold in the case in which the components are not identically distributed, our proof will make essential use of this assumption. We define
	\begin{equation}
		\begin{split}
			R(\Delta(s,t)) &\coloneqq R(t) - R(s) \\
			R(\Delta(s,t),v) &\coloneqq R(t,v) - R(s,v) = \bbE[X_{st} \otimes X_v] \\
			R(\Delta(s,t),\Delta(u,v)) &\coloneqq R(t,v) + R(s,u) - R(t,u) - R(s,v)= \bbE[X_{st} \otimes X_{uv}]
		\end{split}
	\end{equation}
	for $u,v,s,t \in [0,T]$. Note that $R(\Delta(s,t)) \neq R(\Delta(s,t),\Delta(s,t))$.  
	
	We assume $X$ and $R$ satisfy the conditions that make it possible to consider the \emph{signature} of $X$, $\EuS(X)$, defined by the limit in $L^2$ of Stieltjes iterated integrals of smooth or piecewise-linear approximations of $X$, and carry out Malliavin calculus: these are existence of rough path lift and complementary Cameron-Martin regularity \cite[Conditions 2]{CF10} and non-degeneracy of $R$ \cite[Conditions 3]{CF10}. More elementary conditions that imply these may be found, for instance, in \cite{CL19,CL20}. The expected signatures of such processes characterise their law, i.e.\ if $Y$ is any other process with a well-defined signature $\EuS(Y)_{0T}$ (as a $\EuScript{G}(\bbR^d)$-valued random variable) and $\mathbb E \EuS(X)_{0T} = \bbE\EuS(Y)_{0T}$, then $X$ and $Y$ are equal in law: see \cite[Example 6.7]{CL13}, a consequence, among other things, of the greedy estimate \cite{CLL13}. We refer the reader to \cite{CO22} for a treatment of the theory in the case of more general processes, whose expected signatures may not directly characterise the law of the process.
	
	We will denote $\EuS^N(X)$ the signature of $X$ truncated at level $N$ (i.e.\ its projection onto $\bigoplus_{n = 0}^N (\bbR^d)^{\otimes n}$) and $\EuS(X)^{(n)}$ the $n$-th level of the signature (i.e.\ its projection onto $(\bbR^d)^{\otimes n}$). The signature of a process, as that of a path, satisfies two important algebraic relations. The first is the Chen identity, namely that $\EuS(X)_{su} \otimes \EuS(X)_{ut} = \EuS(X)_{st}$. The second is the shuffle identity: letting $\{e_1,\ldots,e_d\}$ denote the canonical basis of $\bbR^d$, and using coordinate notation, i.e.\ $S^{\gamma_1\ldots\gamma_n} \coloneqq \langle e_{\gamma_1} \otimes \cdots \otimes e_{\gamma_n}, S\rangle$ for $S \in T(\!(\bbR^d)\!)$ and $\gamma_1,\ldots,\gamma_n \in [d] \coloneqq \{1,\ldots,d\}$ (and extending linearly), for $0 \leq s \leq t \leq T$ it holds that
	\begin{equation}\label{eq:shuffle}
		\EuS(X)_{st}^{\alpha_1\ldots\alpha_m}\EuS(X)_{st}^{\beta_1\ldots\beta_n} = \EuS(X)_{st}^{(\alpha_1\ldots\alpha_m) \shuffle (\beta_1\ldots\beta_n)}
	\end{equation}
	where $\shuffle$ denotes \say{shuffling} the tuples $\alpha_1\ldots\alpha_m$ and $\beta_1\ldots\beta_n$, i.e.\ summing over all ways of permuting their concatenation $\alpha_1\ldots\alpha_m\beta_1\ldots\beta_n$ whilst preserving the order of each. For further details see, for example, \cite{LCL07}.
	
	In addition to the standard conditions on $R$, we will have to assume a certain amount of smoothness of $R$ together with bounds on its derivatives; the reasons for such hypotheses will be made clear in due course. We assume $R(\,\cdot\,,\,\cdot\,)$ is $C^2$ on the open simplex $\Delta[s,t] \coloneqq \{0 < s < t < T\}$ and continuous on $[0,T]^2$, and that $R(\,\cdot\,)$ is $C^1$ on $(0,T)$. The lack of smoothness assumptions of $R(\,\cdot\,,\,\cdot\,)$ on the diagonal $\{s = t\}$ is crucial for the inclusion of $(1/4,1/2] \ni H$-fBm, which does not even have first partial derivatives on it. Furthermore, we assume there exists an $H \in (0,1)$ with the property that the sample paths of $X$ are either $H$-H\"older, or are $K$-H\"older for all $K < H$; for fBm $H$ will coincide with the Hurst parameter, but the letter $H$ will be used for more general processes to denote the H\"older exponent/supremum of exponents. This the rough path above $X$ will be of finite $1/H$-variation or of finite $p$-variation for all $p> 1/H$.
	
	We also need some quantitative estimates on the derivatives of $R$. Here and throughout the paper, the constant of proportionality implied by the use of $\lesssim$ may only depend on $T, H$ and other general characteristics of the process $X$. We require
	\begin{align}
|\partial_{12}R(s,t)| &\lesssim (t-s)^{2H-2} \qquad \text{on } \Delta[0,T] \label{eq:requiredIneq}\\
\big|\tfrac 12 R'(t) - \partial_2R(s,t)\big| &\lesssim (t-s)^{2H-1} \qquad \text{on } [0,T]^2 \setminus \{s = t\} \label{eq:12REst} \\
|R'(t)| &\lesssim t^{2H-1}  \qquad \qquad \ \ \hspace{0.05em} \text{on } (0,T] \label{eq:R'Est}
\end{align}
	where $\partial_2$ denotes partial differentiation w.r.t.\ the second component and $\partial_{12}$ denotes second-order mixed partial differentiation. Since $R$ is not smooth on the diagonal, the following estimate for on-diagonal square increments of the covariance function, which already appeared in \cite{CQ02}, must be required separately:
	\begin{equation}\label{eq:onDiagEst}
	R(\Delta(s,t), \Delta(s,t)) \lesssim (t-s)^{2H}
	\end{equation} 
	
	We move onto the treatment of Malliavin calculus for $X$. We let $\mathcal H$ be the Hilbert space given by the completion of the following $\bbR$-linear span of \emph{elementary} functions $[0,T] \to \bbR^d$, or equivalently $[0,T] \times [d] \to \bbR$:
	\begin{equation}\label{eq:H}
		\mathcal E \coloneqq \mathrm{span}_\bbR\{ \mathbbm 1^\gamma_{[0,t)} \mid t \in [0,T], \ \gamma = 1,\ldots, d\}
	\end{equation}
	w.r.t.\ the inner product
	\begin{equation}
		\langle \mathbbm 1^\alpha_{[0,s)}, \mathbbm 1^\beta_{[0,t)} \rangle_\mathcal{H} \coloneqq R^{\alpha\beta}(s,t)\text{ .}
	\end{equation}
	Because of independence of components, $\mathcal H$ is equal to an orthogonal direct sum $\mathcal H^1 \oplus \ldots \oplus \mathcal H^d$, and because of equal distribution the direct summands are all equal. Elements of $\mathcal H$ should be viewed as admissible deterministic integrands for $\dif X$, which are represented as Cauchy sequences of elementary integrands in $\mathcal E$. This framework allows us to view the process as an isometry
	\begin{equation}\label{eq:isonormal}
		X \colon \mathcal H \to L^2\Omega,\quad \mathbbm 1^\gamma_{[0,t)} \mapsto X^\gamma_t
	\end{equation}
	often called an isonormal Gaussian process.
	
	The \emph{multiple Wiener integral}
	\begin{equation}
		\skoo^m \colon \mathcal H^{\odot m} \to L^2\Omega 
	\end{equation}
	is the operator defined by the adjoint property (which more generally characterises the divergence operator, when required on random arguments $f$)
	\begin{equation}\label{eq:adjoint}
		\forall Z \in \mathbb D^{m,2} \quad E[ Z \skoo^m(f) ] = \mathbb E[\langle \EuD^m Z, f \rangle_{\mathcal H^{\otimes m}}]
	\end{equation}
	where
	\begin{equation}
		\EuD^m \colon \mathbb D^{m,2} \to L^2(\Omega,\mathcal H^{\odot m})
	\end{equation}
	is the $m^\text{th}$ Malliavin derivative, defined as
	\begin{equation}\label{eq:mal}
		\EuD^m f(X^{\gamma_1}_{t_1}, \ldots, X^{\gamma_n}_{t_n}) \coloneqq \sum_{k_1,\ldots,k_m = 1}^n \partial_{k_1, \ldots ,k_m} f(X^{\gamma_1}_{t_1}, \ldots, X^{\gamma_n}_{t_n}) \mathbbm 1_{[0,t_{k_{\scaleto{1}{3pt}}})}^{\gamma_{k_{\scaleto{1}{3pt}}}} \otimes \cdots \otimes \mathbbm 1_{[0,t_{k_{\scaleto{n}{2pt}}})}^{\gamma_{k_{\scaleto{n}{2pt}}}}
	\end{equation}
	for $f \in C^\infty(\bbR^n)$ with derivatives (including the $0^\text{th}$) of polynomial growth, and extended as a closed operator to a certain domain $\mathbb D^{m,2}$. $\mathcal H^{\odot m}$ denotes the subspace of $\mathcal H^{\otimes m}$ (the tensor product taken in the category of Hilbert spaces) of symmetric tensors. $\EuD^m$ takes a square-integrable random variable and returns a random element of $\mathcal H^{\odot m}$, which in case of membership to $\mathcal E^{\odot m}$ (or otherwise a function member of $\mathcal H^{\odot m}$ in the sense of \autoref{def:fH} below) will be a function of $m$ $($time$,$index$)$ pairs. Note that, while $\skoo$ is symmetric in the sense that it is left invariant by permuting $($time$,$index$)$ jointly, it is not symmetric if only time variables or indices are permuted (e.g.\ it is possible to use $\skoo$ to define a Lévy area --- see \autoref{expl:w3} below). When $\EuD^mZ$ is a function, as in the case \eqref{eq:mal}, we denote its evaluation on $m$ $($time$,$index$)$ pairs $\EuD_{(u_1,\gamma_1),\ldots,(u_m,\gamma_m)}Z$; occasionally it may make more sense to suppress the indices in the notation, in which case we can just write  $\EuD_{u_1,\ldots,u_m}Z$. We may extend $\skoo$ to a map $\skoo^m \coloneqq \mathcal H^{\otimes m} \to L^2\Omega$ by pre-composing with symmetrisation, and we have for $f,g \in \mathcal H^{\otimes m}$
	\begin{equation}
		\mathbb E[\skoo^m(f) \skoo^n(g)]=  \kron^{mn} \sum_{\sigma \in \mathfrak S_m} \langle f, \sigma_*g \rangle_{\mathcal H^{\otimes m}}\text{ .}
	\end{equation}
	This implies that multiple Wiener integration defines an isometry
	\begin{equation}\label{eq:skoIsometry}
		\skoo^{\scriptscriptstyle \bullet} \colon \bigoplus_{m = 0}^\infty \mathcal H^{\odot m} \xrightarrow{\cong} L^2 \Omega
	\end{equation}
	where the source is given the degree-wise rescaled inner product $(f,g) \mapsto m!\langle f,g \rangle_{\mathcal H^{\otimes m}}$ for $f,g$ of the same degree and zero otherwise, and $\Omega$ is endowed with the sigma algebra generated by the process $X_{t \in [0,T]}$. The image of the $m$-th Wiener integral operator, the space of the random variables $\skoo^m(f)$ with $f$ ranging in $\mathcal H^{\odot m}$, is called the $m$-th \emph{Wiener chaos} of $X$. We denote it $\mathscr W^m$ and the $m$-th Wiener chaos projection $\mathscr w^m \colon L^2 \Omega \twoheadrightarrow \mathscr W^m$. Note that $\mathscr w^0 = \mathbb E$ with values in $\mathscr W^0 = \bbR$, while $\mathscr W^1$ is given by linear functionals of $X$. We thus have the Wiener chaos decomposition $L^2 \Omega = \bigoplus_{m = 0}^\infty \mathscr W^m$
	which means it is possible to represent any random variable in $L^2\Omega$ (measurable w.r.t.\ to the sigma-algebra generated by $X$) as an $L^2$-absolutely convergent series
	\begin{equation}
		L^2\Omega \ni Z = \sum_{m = 0}^\infty \mathscr w^m Z,\qquad \lVert Z \rVert_{L^2}^2 = \sum_{m = 0}^\infty \lVert \mathscr w^m Z \rVert_{L^2}^2 = \sum_{m = 0}^\infty m! \lVert f^m \rVert_{\mathcal H^{\otimes m}}^2
	\end{equation}
	where $f^m = (\skoo^m)^{-1} \circ \mathscr w^m (Z)$. The map $(\skoo^m)^{-1} \circ \mathscr w^m$ admits an expression in terms of the Malliavin derivative: this is \emph{Stroock's formula}, which states that for $Z \in \mathbb D^{m,2}$
	\begin{equation}\label{eq:stroock}
		(\skoo^m)^{-1} \circ \mathscr w^m (Z) = \frac{1}{m!} \mathbb E[\EuD^m Z] \text{ .}
	\end{equation}
	As a consequence, if $Z \in \mathbb D^{\infty,2} \coloneqq \bigcap_{m = 0}^\infty \mathbb D^{m,2}$ we can write its Wiener chaos decomposition as the series
	\begin{equation}
		Z = \sum_{m = 0}^\infty \frac{1}{m!} \skoo^m \mathbb E[ \EuD^m Z]\text{ .}
	\end{equation}

 We continue calling elements of $\mathcal E^{\otimes m}$ elementary functions, in light of the fact that they can be identified with functions $([0,T] \times [d])^m \to \bbR$ by the mapping
	\begin{equation}\label{eq:idfun}
		\mathbbm 1^{\gamma_1}_{[0,t_1)} \otimes \cdots \otimes \mathbbm 1^{\gamma_m}_{[0,t_m)} \mapsto \mathbbm 1^{\gamma_1,\ldots,\gamma_m}_{[0,t_1) \times \cdots \times [0,t_m)}\text{ .}
	\end{equation}
	This is the map given by the product of the Kronecker deltas $\kron^{\gamma_1}_\cdot \cdots \kron^{\gamma_m}_\cdot$ and the indicator function on the $m$-cube $[0,t_1) \times \cdots \times [0,t_m)$, each $\kron$ paired with the respective time variable. Since $\mathcal E^{\otimes m}$ is dense in $\mathcal H^{\otimes m}$, elements of the latter may be identified as equivalence classes of Cauchy sequences in $\mathcal E^{\otimes m}$. While $\mathcal H^{\otimes m}$ is not, in general, a space of functions, it is possible to uniquely associate elements of $\mathcal H^{\otimes m}$ to certain measurable functions $([0,T] \times [d])^m \to \bbR$ as follows:
	\begin{defn}[Functions as elements of $\mathcal H^{\otimes m}$]\label{def:fH}
	For a function $f \colon ([0,T] \times [d])^m \to \bbR$ we will write $f \in \mathcal H^{\otimes m}$ if there exist a Cauchy sequence $(f_n)_n \subset \mathcal E^{\otimes m}$, uniformly bounded as a sequence of functions (according to the identification \eqref{eq:idfun}), with $f_n \to f$ a.e. In this case we will say that $f$ \emph{represents} $\lim f_n \in \mathcal H^{\otimes m}$. If $f$ represents $\phi, \psi \in \mathcal H^{\otimes m}$ then $\phi = \psi$: this is an immediate consequence of the following
	\end{defn}
\begin{lem}
Let $(f_n)_n$ be as in the above definition with $f = 0$. Then $f_n \to 0$ in $\mathcal H^{\otimes m}$.
\begin{proof}
	Let 
	\[
	f_n = \sum_{\gamma_1,\ldots,\gamma_m = 1}^d f_{n;\gamma_1,\ldots,\gamma_m} \mathbbm 1^{\gamma_1,\ldots,\gamma_m}
	\]
	with $f_{n;\gamma_1,\ldots,\gamma_m} \colon [0,T]^m \to \bbR$. Then $f_n \to 0$ a.e.\ if and only if $f_{n;\gamma_1,\ldots,\gamma_m} \to 0$ a.e.\ for each $(\gamma_1,\ldots,\gamma_m) \in [d]^m$. Keeping in mind that $\mathcal H \cong (\mathcal H^1)^{\bigoplus d}$ we may therefore assume $d = 1$ and suppress indices. Following \cite[p.588]{HuCa78}, we test the sequence with elementary functions: letting $f_n = \sum_{s^n_1,\ldots,s^n_m} f_n^{s^n_1,\ldots,s^n_m}\mathbbm 1_{[0,s^n_1)\times \ldots\times[0,s^n_m)}$ and $\mathcal E^{\otimes m} \ni g = \sum_{t_1,\ldots,t_m} g^{t_1,\ldots,t_m}\mathbbm 1_{[0,t_1)\times\ldots\times[0,t_m)}$ with $f_n^{s^n_1,\ldots,s^n_m}, g^{t_1,\ldots,t_m} \in \bbR$ uniformly bounded (and the sums finite) we have that
	\begin{align*}
		&\langle f_n, g \rangle_{\mathcal H^{\otimes m}} \\
		={} &\sum_{\substack{s^n_1,\ldots,s^n_m \\ t_1,\ldots,t_m}}f_n^{s^n_1,\ldots,s^n_m} g^{t_1,\ldots,t_m} R(s^n_1,t_1) \cdots R(s^n_m,t_m)  \\
		 ={} &\sum_{t_1,\ldots,t_m} g^{t_1,\ldots,t_m} \int_{((0,T] \setminus \{t_1\}) \times \cdots \times ((0,T] \setminus \{t_m\})}  f_n(s_1,\ldots,s_m) \partial_1 R(s_1,t_1) \cdots \partial_1 R(s_m,t_m) \dif s_1 \cdots \dif s_m\text{ .}
	\end{align*}
\eqref{eq:12REst} and \eqref{eq:R'Est} imply that the integrands are absolutely and uniformly bounded by $ [|t_1-s_1|^{2H-1} \vee s_1^{2H-1}] \cdots [|t_m-s_m|^{2H-1} \vee s_m^{2H-1}]$ (up to a constant), which is integrable on $((0,T] \setminus \{t_1\}) \times \cdots \times ((0,T] \setminus \{t_m\})$. By dominated convergence $\langle \phi, g \rangle_{\mathcal H^{\otimes m}} = \lim\langle f_n, g \rangle_{\mathcal H^{\otimes m}} = 0$, where $\phi \coloneqq \lim f_n$ in $\mathcal H^{\otimes m}$, and $\phi = 0$ follows from the fact that $g$ ranges in a dense set. 
\end{proof}
\end{lem}
	
In light of the aforementioned non-degeneracy condition on $X$, we also expect the converse to hold: if $\phi \in \mathcal H^{\otimes m}$ is represented by the functions $f,g$ in the above sense, then $f = g$ a.e. An example of a degenerate stochastic process, for which this property would not hold, is given by taking any process $X$ and concatenating it with itself path by path; the resulting covariance function $R$ would be invariant under transposing the intervals $[0,T)$ and $[T,2T)$. We also note that, in specific cases, it is possible to describe $\mathcal H$ explicitly: if $X$ is a fractional Brownian motion with Hurst parameter $H \in (0,1)$, the identity on $\mathcal E$ induces an isomorphism between $\mathcal H$ and the Sobolev space $W^{1/2-H, 2}$ \cite{Jol07}, which is a space of functions for $H \in (0,1/2]$ but not for $H \in (1/2,1)$.

We will mostly be considering Wiener integrals on simplices, which has the effect of quotienting out symmetry of the operator $\skoo^m$. We will often resort to integral notation, e.g.\ if $\mathbbm 1^{\alpha\beta}_{\Delta[s,t]} \in \mathcal H^{\otimes 2}$ (the function that maps $((u,\gamma),(v,\delta)) \mapsto \kron^{\alpha\gamma}\kron^{\beta\delta} \mathbbm 1_{s < u < v < t}$) in the sense of \autoref{def:fH}, we will write $\skoo^2(\mathbbm 1^{\alpha\beta}_{\Delta[s,t]}) \eqqcolon \int_{s < u < v < t} \sko X_u^\alpha \sko X^\beta_v$ to be the limit in $L^2$ of $\skoo^2(f_n)$. Wiener integrals of elements of $\mathcal E^{\otimes m}$, on the other hand, can be computed explicitly by using the adjoint property \eqref{eq:adjoint}: for example, it can be checked that
	\[
	\skoo^2(\mathbbm 1^{\alpha}_{[s,t)} \otimes  \mathbbm 1^{\beta}_{[u,v)}) = X^{\alpha}_{st}  X^{\beta}_{uv} - R^{\alpha\beta}(\Delta(s,t),\Delta(u,v))\text{ .}
	\]
	The more general formula involves multivariate analogues of the Hermite polynomials (see \cite[\S 2.7.2]{NouPec12} and \cite[p.244]{Fer22}). When $X$ is a Gaussian martingale (but not necessarily if it is only a semimartingale), multiple Wiener integration on the simplex coincides with iterated Wiener-It\^o integration.

	\section{The main result, some consequences}\label{sec:main}
	
	We begin this section with some more notation. We denote $[n] \coloneqq \{1,\ldots,n\}$ the set with $n$ elements. We will be concerned with iterated integrals on the $n$-simplex $\Delta^n[s,t] \coloneqq \{(u_1,\ldots,u_n) \mid s < u_1 < \ldots < u_n < t \}$. Because such integrals will involve the covariance function, integration variables will sometimes come in pairs. For $m, n \in \mathbb N$ we denote $\mathcal P^n_m$ the collection of partitions of subsets of $[n]$ of cardinality $n-m$ into sets of cardinality $2$. Note that this means $\mathcal P^n_m = \varnothing$ whenever $n \neq m \ (\mathrm{mod} \ 2)$ or $m > n$, but $\mathcal P^n_n$ has precisely one element, $\varnothing$: the empty set admits the empty collection of subsets as a partition, which vacuously belongs to $\mathcal P^n_n$. For example, $Q \coloneqq \{\{1,4\},\{3,8\},\{5,6\}\} \in \mathcal P^8_2$ viewed as a partition of the set $\{1,3,4,5,6,8\} \subseteq [8]$. For $P \in \mathcal P^n_m$ we will denote $\overline P \coloneqq [n] \setminus \cup P$ (in the partition of the above example, $\overline Q = \{2,7\}$).
	
	It will convenient to use graphical notation to denote such objects, and for reasons that will become apparent shortly, for a pair $\{i,j\}$ with $i \leq j$ we will distinguish between the consecutive case $j = i+1$ and the non-consecutive one $j > i+1$. The partition $Q \in \mathcal P^8_2$ above is represented by
	\begin{equation}\label{eq:diagEx}
		Q = \tikz[baseline = -0.1ex]{
			\draw[fill] (0,0) circle [radius=0.04];
			\draw[fill] (0.3,0) circle [radius=0.04];
			\draw[fill] (0.6,0) circle [radius=0.04];
			\draw[fill] (0.9,0) circle [radius=0.04];
			\draw[fill] (1.2,0) circle [radius=0.04];
			\draw[fill] (1.5,0) circle [radius=0.04];
			\draw[fill] (1.8,0) circle [radius=0.04];
			\draw (0,0) .. controls (0,0.4) and (0.9,0.4) .. (0.9,0);
			\draw (0.6,0) .. controls (0.6,0.4) and (1.8,0.4) .. (1.8,0);
			\draw (1.2,0)-- (1.2,0.2);
		}\quad \text{.}
	\end{equation}
	We will refer to such graphics as \emph{diagrams}. We have drawn one node for each $i \in [n]$ that is not paired with a consecutive integer, and one node for each consecutive pair (in this case only $\{5,6\}$); when counting nodes, a node corresponding to such a pair should be thought as having double weight. In our example, the $5^\text{th}$ node actually counts for positions $5$ and $6$. With this convention, for each non-consecutive pair $\{i,j\}$ we have drawn an \emph{arc} connecting the two nodes of positions $i$ and $j$, and for each node corresponding to a consecutive pair we have drawn a \emph{line} going upwards. Nodes that do not have a line or arc entering them correspond to elements of $\overline P$, and we will call them \emph{single}. Note that, by construction, there is never an arc between two consecutive nodes: this will be critical for convergence of the associated integrals described below. In the next section, we will be particularly concerned with maximal \emph{sequences of consecutive pairings}, i.e.\ collections of pairings $\{k,k+1\},\ldots,\{k+l,k+l+1\} \in P$ with $l \geq 0$ and s.t.\ $\{k-2,k-1\},\{k+l+2,k+l+3\} \not \in P$.
	
	Now, given $P \in \mathcal P^n_m$, $0 \leq s \leq t \leq T$ and $\gamma_1,\ldots,\gamma_n \in [d]$ we associate to it a continuous function $P_{st}^{\gamma_1,\ldots,\gamma_n} \colon \Delta^m[s,t] \times [d]^m \to \bbR$ by integrating over as many variables as there are non-single nodes in the diagram that represents $P$: call this number, which equals twice the number of non-consecutive pairs in $P$ plus the number of consecutive ones, $\#P$. This explains our choice for the above notation: each node either corresponds to an integration variable or to a free variable, i.e.\ a variable of which $P_{st}^{\gamma_1,\ldots,\gamma_n}$ is a function. We use the shorthands
	\begin{equation}\label{eq:difCases}
		\begin{split}
			R(\dif u_i, \dif u_j) &\coloneqq \partial_{12}R(u_i,u_j) \dif u_i \dif u_j  \\ \tfrac 12 R(\dif u_{h+1}) - R(u_{h-1},\dif u_{h+1}) &\coloneqq \big[\tfrac 12 R'(u_{h+1}) - \partial_{2}R(u_{h-1},u_{h+1})\big] \dif u_{h+1}\text{ ,}
		\end{split}
	\end{equation}
	and the former will only be used when $j > i+1$. Crucially, we are defining the second case as $\frac 12 R(\dif u_{h+1}) - R(u_{h-1},\dif u_{h+1})$, not as $R(u_{h+1},\dif u_{h+1}) - R(u_{h-1},\dif u_{h+1})$, since this would be ill-defined in many cases (including $1/2>H$-fBm) because $R(\,\cdot\,, \,\cdot\,)$ may not admit partial derivatives on the diagonal. On the other hand, we are assuming that the variance function $R(\,\cdot\,)$ is differentiable. 
	\begin{defn}[$P^{\gamma_1,\ldots,\gamma_n}_{st}$]\label{def:defIntegral}
		For $\gamma_1,\ldots,\gamma_n \in [d]$, $0 \leq s \leq t \leq T$ and $P \in \mathcal P^n_m$ define
		\begin{equation}\label{eq:defIntegral}
			\begin{split}
				P_{st}^{\gamma_1,\ldots,\gamma_n}(u_k \mid k \in \overline P)
				&\coloneqq  \prod_{k \in \overline P} \mathbbm 1^{\gamma_k} \cdot \int_{\Delta^{\#P}[s,t]} \prod_{\substack{\{i,j\} \in P \\ |j-i|>1}} R^{\gamma_i \gamma_j}(\dif u_i,\dif u_j)  \\ &\qquad \quad \qquad \quad \cdot  \prod_{\{h,h+1\} \in P} \big[ \tfrac 12 R^{\gamma_h \gamma_{h+1}}(\dif u_{h+1}) - R^{\gamma_h \gamma_{h+1}}(u_{h-1},\dif u_{h+1}) \big]
			\end{split}
		\end{equation}
	as a function $([0,T] \times [d])^m \to \bbR$ extended with the value $0$ outside $\Delta^m[s,t]$.
	\end{defn}
	The variables $u_k$ with $k \in \overline P$ are supplied as arguments, so in fact this is an integral over a disjoint union of up to $m+1$ simplices (fewer if some of the elements of $\overline P$ are consecutive). The $k^\text{th}$ index in $[d]^n$ is given as argument to $\mathbbm 1^{\gamma_k}$ as a Kronecker delta: this means that $P_{st}^{\gamma_1,\ldots,\gamma_n}$ vanishes on all but one element of $[d]^m$. The reason why we still consider $P_{st}^{\gamma_1,\ldots,\gamma_n}$ as a function on $[d]^m$ is that this is necessary to view it as an element of $\mathcal H^{\otimes m}$; nevertheless, when the indices are fixed it will sometimes be convenient to just think of it as a function of $m$ times. If $m = 0$, $P^{\gamma_1,\ldots,\gamma_n}_{st}$ is just a real number.
	\begin{rem}
		The presence of the second type of integrand in \autoref{def:defIntegral} is the reason for the smoothness assumptions on the variance and covariance functions, which are not to be found in most of the literature on these topics: this is because it would be difficult to define integrals such as $\int_{s < u < v < t} \big[\tfrac 12 R(\dif u) -  R(s,\dif u)\big]\big[\tfrac 12 R(\dif v) - R(u,\dif v)\big]$ as iterated Young integrals, without taking derivatives, since the variable $u$ in its undifferentiated form appears after the integrator $\frac 12 R(\dif u) -  R(s,\dif u)$; this of course is no longer an issue under our smoothness hypotheses, thanks to which the above integral is defined as the Lebesgue integral on the simplex $\int_{s < u < v < t} \big[\tfrac 12 R'(u) - \partial_2 R(s,u)\big]\big[\tfrac 12 R'(v) - \partial_2 R(u,v)\big] \dif u \dif v$.
	\end{rem}
	
	When $P$ is represented by a diagram, we will decorate the nodes with labels. For example, the integral associated to \eqref{eq:diagEx} with labelling $\alpha,\ldots,\vartheta$ is given by
	\begin{equation*}
		\begin{split}
			&(\!\tikz[baseline = -0.1ex]{
				\draw[fill] (0,0) circle [radius=0.04];
				\draw[fill] (0.3,0) circle [radius=0.04];
				\draw[fill] (0.6,0) circle [radius=0.04];
				\draw[fill] (0.9,0) circle [radius=0.04];
				\draw[fill] (1.2,0) circle [radius=0.04];
				\draw[fill] (1.5,0) circle [radius=0.04];
				\draw[fill] (1.8,0) circle [radius=0.04];
				\draw (0,0) .. controls (0,0.4) and (0.9,0.4) .. (0.9,0);
				\draw (0.6,0) .. controls (0.6,0.4) and (1.8,0.4) .. (1.8,0);
				\draw (1.2,0)-- (1.2,0.2);
				\node[below] at (0,0) {$\scriptstyle \alpha$};
				\node[below] at (0.3,0) {$\scriptstyle \beta$};
				\node[below] at (0.6,0) {$\scriptstyle \gamma$};
				\node[below] at (0.9,0) {$\scriptstyle \delta$};
				\node[below] at (1.2,0) {$\scriptstyle \varepsilon \hspace{-0.025em} \zeta$};
				\node[below] at (1.5,0) {$\scriptstyle  \eta$};
				\node[below] at (1.8,0) {$\scriptstyle \vartheta$};
			}\!)_{st} \\
			={} & \mathbbm 1^{\beta\eta} \int_{\Delta^7[s,t]} R^{\alpha\delta}(\dif u_1,\dif u_4)R^{\gamma\vartheta}(\dif u_3, \dif u_8)\big[\tfrac 12 R^{\varepsilon\zeta}(\dif u_6) - R^{\varepsilon\zeta}(u_4, \dif u_6)\big] \\
			={} & \mathbbm 1_{\Delta^2[s,t]}^{\beta\eta}(u_2,u_6)  \int_{\substack{s < u_1 < u_2 \\ u_2 < u_3 < u_4 < u_6 < u_7 \\ u_7 < u_8 < t}} \partial_{12}R^{\alpha\delta}(u_1,u_4)\partial_{12}R^{\gamma\vartheta}(u_3,u_8)\big[\tfrac 12 R^{\varepsilon\zeta}{}'(u_6) - \partial_2R^{\varepsilon\zeta}(u_4, u_6)\big] \\[-1.5em]
			&\hspace{16em} \dif u_1 \dif u_3 \dif u_4 \dif u_6 \dif u_8 \text{ .}
		\end{split}
	\end{equation*}
	This is viewed as a function of the variables $u_2,u_6$ ranging on the simplex $\Delta^2[s,t]$, each paired with an index variable, which must respectively be equal to $\beta$, $\eta$ for the expression not to vanish. The variable $u_5$ has been skipped, since it is the first term in the consecutive pair $\{5,6\}$. We will show that integrals defined in this fashion are a.e.\ limits of Cauchy sequences in $\mathcal E^{\otimes m}$, which therefore uniquely represent elements of $\mathcal H^{\otimes m}$ according to \autoref{def:fH}. When taking multiple Wiener integrals of them, the indices corresponding to the nodes that represent free variables will become the coordinate processes that are being integrated against, e.g.
	\begin{equation*}
		\begin{split}
			&\skoo^2 (\!\tikz[baseline = -0.1ex]{
				\draw[fill] (0,0) circle [radius=0.04];
				\draw[fill] (0.3,0) circle [radius=0.04];
				\draw[fill] (0.6,0) circle [radius=0.04];
				\draw[fill] (0.9,0) circle [radius=0.04];
				\draw[fill] (1.2,0) circle [radius=0.04];
				\draw[fill] (1.5,0) circle [radius=0.04];
				\draw[fill] (1.8,0) circle [radius=0.04];
				\draw (0,0) .. controls (0,0.4) and (0.9,0.4) .. (0.9,0);
				\draw (0.6,0) .. controls (0.6,0.4) and (1.8,0.4) .. (1.8,0);
				\draw (1.2,0)-- (1.2,0.2);
				\node[below] at (0,0) {$\scriptstyle \alpha$};
				\node[below] at (0.3,0) {$\scriptstyle \beta$};
				\node[below] at (0.6,0) {$\scriptstyle \gamma$};
				\node[below] at (0.9,0) {$\scriptstyle \delta$};
				\node[below] at (1.2,0) {$\scriptstyle \varepsilon \hspace{-0.025em} \zeta$};
				\node[below] at (1.5,0) {$\scriptstyle  \eta$};
				\node[below] at (1.8,0) {$\scriptstyle \vartheta$};
			}\!)_{st} \\
			={}&\int_{s < u_2 < u_7 < t} \Bigg[ \int_{\substack{s < u_1 < u_2 \\ u_2 < u_3 < u_4 < u_6 < u_7 \\ u_7 < u_8 < t}} \partial_{12}R^{\alpha\delta}(u_1,u_4)\partial_{12}R^{\gamma\vartheta}(u_3,u_8)\big[\tfrac 12 R'^{\varepsilon\zeta}(u_6) - \partial_2R^{\varepsilon\zeta}(u_4, u_6)\big] \Bigg] \sko X^\beta_{u_2} \sko X^\eta_{u_7}\text{ .} \\[-1.5em]
			&\hspace{14em} \dif u_1 \dif u_3 \dif u_4 \dif u_6
		\end{split}
	\end{equation*}
	We are now ready to state the main theorem.
	
	\begin{thm}[Wiener chaos expansion of the signature of a Gaussian process]\label{thm:main} \ \\
		Given $m,n \in \mathbb N$, $P \in \mathcal P^n_m$, $\gamma_1,\ldots,\gamma_n \in [d]$, $0 \leq s \leq t \leq T$, it holds that $P^{\gamma_1,\ldots,\gamma_n}_{st} \in \mathcal H^{\otimes m}$ in the sense of \autoref{def:fH}, and the $m^\text{th}$ Wiener chaos projection of the signature of $X$ is given by
		\begin{equation}
			\mathscr w^m \EuS(X)^{\gamma_1,\ldots,\gamma_n}_{st} =  \sum_{P \in \mathcal P^n_m} \skoo^m P^{\gamma_1,\ldots,\gamma_n}_{st}\text{ .}
		\end{equation}
	\end{thm}
	In particular, notice that $\mathscr w^m \EuS(X)^{\gamma_1,\ldots,\gamma_n}_{st}$ can only be non-zero when $m \leq n$ and $m \equiv n \ (\mathrm{ mod \ 2})$. The most important case of this result is when $m = 0$: 
	\begin{cor}[Expected signature of a Gaussian process]\label{cor:esig}
		With notation as above, we have
		\begin{equation}\label{eq:esig}
			\mathbb E \EuS(X)^{\gamma_1,\ldots,\gamma_n}_{st} = \sum_{P \in \mathcal P^n_0} P^{\gamma_1,\ldots,\gamma_n}_{st}\text{ .}
		\end{equation}
	\end{cor}
	
	\begin{rem}[Eliminating variables]\label{rem:simpl}
		While convergence rules out always considering integrands of the first type in \eqref{eq:difCases} (which would mean allowing diagrams with arcs between consecutive nodes), one may wonder whether it is possible to only consider integrands of the second type, i.e.\ by integrating out one variable per pair and thus simplifying the presentation of the formula. This, however, is not possible in general, because of the additional constraint that requires two consecutive variables not to be both integrated out (for the expression to make sense as an integral). It is not difficult to see, for example, that in the following diagram
		\[
		\tikz[baseline = -0.1ex]{
			\draw[fill] (0,0) circle [radius=0.04];
			\draw[fill] (0.3,0) circle [radius=0.04];
			\draw[fill] (0.6,0) circle [radius=0.04];
			\draw[fill] (0.9,0) circle [radius=0.04];
			\draw[fill] (1.2,0) circle [radius=0.04];
			\draw[fill] (1.5,0) circle [radius=0.04];
			\draw (0,0) .. controls (0,0.5) and (1.5,0.5) .. (1.5,0);
			\draw (0.3,0) .. controls (0.3,0.3) and (0.9,0.3) .. (0.9,0);
			\draw (0.6,0) .. controls (0.6,0.3) and (1.2,0.3) .. (1.2,0);
		}
		\]
		at most two variables can be integrated out (unless the remaining integral can be solved or simplified analytically). Luckily, the only case in which it is necessary for convergence to integrate out certain variables (as specified in the second case of \eqref{eq:difCases}), is when there are consecutive pairs: this is always possible, even when more than one pair in a row is consecutive, since we may always pick the first variable to integrate out (as done here --- one could equivalently have chosen the second). Of course, there is always some number of additional variables that can be eliminated, but we do not immediately see a way of doing this in a maximal way that is canonical.
	\end{rem}

	\begin{expl}[The Wiener chaos decomposition of $\EuS^3(X)_{st}$]\label{expl:w3}
		We give the explicit expression for the Wiener chaos expansion of the signature truncated at level $3$. These terms are especially significant, considering that they are the ones that define the rough path when $1/4 < H \leq 1/3$: higher signature terms can be derived in a pathwise fashion by Lyons's extension theorem without involving probability. We represent each signature term as a sum of their Wiener chaos projections in ascending order; in particular the sum of all non-random terms constitutes the expectation of the left hand side.
		\begin{align*}
			\EuS(X)^{\varnothing}_{st} &= \varnothing_{st} = 1, \quad 	\EuS(X)^{\gamma}_{st} = \skoo(\!\tikz[baseline = -0.1ex]{
				\draw[fill] (0,0) circle [radius=0.04];
				\node[below] at (0,0) {$\scriptstyle \gamma$};
			}\!)_{st} = \skoo (\mathbbm 1_{[s,t)}^\gamma) = X_{st}^\gamma \\[1em]
			\EuS(X)^{\alpha\beta}_{st} &= \tikz[baseline = -0.1ex]{
				\draw[fill] (0,0) circle [radius=0.04];
				\node[below] at (0,0) {$\scriptstyle \alpha \hspace{-0.05em} \scriptstyle \beta$};
				\draw (0,0) -- (0,0.3);
			} + \skoo^{2}(\!\tikz[baseline = -0.1ex]{
				\draw[fill] (0,0) circle [radius=0.04];
				\draw[fill] (0.3,0) circle [radius=0.04];
				\node[below] at (0,0) {$\scriptstyle \alpha$};
				\node[below] at (0.3,0) {$\scriptstyle \beta$};
				\draw[fill] (0,0) circle [radius=0.04];
			} \!) \\
			&= \frac{R^{\alpha\beta}(s) + R^{\alpha\beta}(t)}{2} - R^{\alpha\beta}(s,t) + \int_{s < u < v < t} \sko X^\alpha_u \sko X^\beta_v \\[1em]
			\EuS(X)^{\alpha\beta\gamma}_{st} &=  \skoo(\!\tikz[baseline = -0.1ex]{
				\draw[fill] (0,0) circle [radius=0.04];
				\draw[fill] (0.3,0) circle [radius=0.04];
				\node[below] at (0,0) {$\scriptstyle \alpha\hspace{-0.05em}\scriptstyle \beta$};
				\node[below] at (0.3,0) {$ \scriptstyle \gamma$};
				\draw (0,0) -- (0,0.3);
				\draw[fill] (0,0) circle [radius=0.04];
			} \!)_{st} + \skoo(\!\tikz[baseline = -0.1ex]{
				\draw[fill] (0,0) circle [radius=0.04];
				\draw[fill] (0.3,0) circle [radius=0.04];
				\node[below] at (0,0) {$\scriptstyle \alpha$};
				\node[below] at (0.3,0) {$\scriptstyle \beta\hspace{-0.05em} \scriptstyle \gamma$};
				\draw (0.3,0) -- (0.3,0.3);
				\draw[fill] (0,0) circle [radius=0.04];
			} \!)_{st} + \skoo(\!\tikz[baseline = -0.1ex]{
				\draw[fill] (0,0) circle [radius=0.04];
				\draw[fill] (0.3,0) circle [radius=0.04];
				\draw[fill] (0.6,0) circle [radius=0.04];
				\node[below] at (0,0) {$\scriptstyle \alpha$};
				\node[below] at (0.3,0) {$\scriptstyle \beta$};
				\node[below] at (0.6,0) {$\scriptstyle \gamma$};
				\draw[fill] (0,0) circle [radius=0.04];
				\draw (0,0) .. controls (0,0.3) and (0.6,0.3) .. (0.6,0);
			} \!)_{st} + \skoo^{3}(\!\tikz[baseline = -0.1ex]{
				\draw[fill] (0,0) circle [radius=0.04];
				\draw[fill] (0.3,0) circle [radius=0.04];
				\draw[fill] (0.6,0) circle [radius=0.04];
				\node[below] at (0,0) {$\scriptstyle \alpha$};
				\node[below] at (0.3,0) {$\scriptstyle \beta$};
				\node[below] at (0.6,0) {$\scriptstyle \gamma$};
				\draw[fill] (0,0) circle [radius=0.04];
			} \!)_{st} \\
			&= \int_s^t \bigg( \frac{R^{\alpha\beta}(s) + R^{\alpha\beta}(u)}{2} - R^{\alpha\beta}(s,u) \bigg)\sko X^\gamma_u +  \int_s^t \bigg( \frac{R^{\beta\gamma}(u) + R^{\beta\gamma}(t)}{2} - R^{\beta\gamma}(u,t) \bigg)\sko X^\alpha_u\\
			&\phantom{{}={}} + \int_s^t R^{\alpha\gamma} (\Delta(s,u),\Delta(u,t)) \sko X^\beta_u +  \int_{s < u < v < w <t} \sko X^\alpha_u \sko X^\beta_v \sko X^\gamma_w
		\end{align*}
		In particular, notice how the expected signature of level $2$ is given by the difference between the average of the variances and the covariance:
		\begin{equation}\label{eq:ES2}
			\mathbb E \EuS(X)^{\alpha\beta}_{st} = \frac{R^{\alpha\beta}(s) + R^{\alpha\beta}(t)}{2} - R^{\alpha\beta}(s,t)
		\end{equation}
		and that the statement that \say{the It\^o and Stratonovich Lévy areas are equal} carries over to the Gaussian Wiener-rough setting, in the sense that
		\begin{equation}\label{eq:Levy}
			\frac 12 \big(\EuS(X)^{\alpha\beta}_{st}-\EuS(X)^{\beta\alpha}_{st} \big) =  \frac 12 \int_{s < u < v < t} \sko X^\alpha_u \sko X^\beta_v - \sko X^\beta_u \sko X^\alpha_v 
		\end{equation}
		by symmetry of the covariance function.
	\end{expl}

	\begin{expl}[$\bbE\EuS(X)^{(4)}$]\label{expl:EX4}
		\autoref{cor:esig} at level $4$ is given by
		\begin{align}\label{eq:level4integrals}
			\begin{split}
				\bbE \EuS(X)^{\alpha\beta\gamma\delta}_{st} ={}&(\!\tikz[baseline = -0.1ex]{
					\draw[fill] (0,0) circle [radius=0.04];
					\draw[fill] (0.3,0) circle [radius=0.04];
					\node[below] at (0,0) {$\scriptstyle \alpha\hspace{-0.05em}\scriptstyle \beta \ $};
					\node[below] at (0.3,0) {$\ \scriptstyle \gamma\hspace{-0.05em}\scriptstyle \delta$};
					\draw (0,0)-- (0,0.3);
					\draw (0.3,0)-- (0.3,0.3);
				}\! )_{st} + (\!\tikz[baseline = -0.1ex]{
					\draw[fill] (0,0) circle [radius=0.04];
					\draw[fill] (0.3,0) circle [radius=0.04];
					\draw[fill] (0.6,0) circle [radius=0.04];
					\node[below] at (0,0) {$\scriptstyle \alpha$};
					\node[below] at (0.3,0) {$\scriptstyle \beta\hspace{-0.05em}\scriptstyle \gamma$};
					\node[below] at (0.6,0) {$\scriptstyle \delta$};
					\draw (0,0) .. controls (0,0.4) and (0.6,0.4) .. (0.6,0);
					\draw (0.3,0)-- (0.3,0.2);
				} \!)_{st} +  (\! \tikz[baseline = -0.1ex]{
					\draw[fill] (0,0) circle [radius=0.04];
					\draw[fill] (0.3,0) circle [radius=0.04];
					\draw[fill] (0.6,0) circle [radius=0.04];
					\draw[fill] (0.9,0) circle [radius=0.04];
					\node[below] at (0,0) {$\scriptstyle \alpha$};
					\node[below] at (0.3,0) {$\scriptstyle \beta$};
					\node[below] at (0.6,0) {$\scriptstyle \gamma$};
					\node[below] at (0.9,0) {$\scriptstyle \delta$};
					\draw (0,0) .. controls (0,0.3) and (0.6,0.3) .. (0.6,0);
					\draw (0.3,0) .. controls (0.3,0.3) and (0.9,0.3) .. (0.9,0);}\!)_{st} \\
				={}&\int_{s < u < v < t} \big[ \tfrac 12 R^{\alpha\beta}(\dif u) - R^{\alpha\beta}(s,\dif u) \big] \big[ \tfrac 12 R^{\gamma\delta}(\dif v) - R^{\gamma\delta}(u,\dif v) \big]\\
				&+\int_{s < u < v < w < t} R^{\alpha\delta}(\dif u,\dif w) \big[ \tfrac 12 R^{\beta\gamma}(\dif v) - R^{\beta\gamma}(u,\dif v) \big] \\
				&+\int_{s < u < v < w < z < t}  R^{\alpha\gamma}( \dif u, \dif w)  R^{\beta\delta}( \dif v, \dif z)\text{ .}
			\end{split}
		\end{align}
		Using a clever transformation, \cite[Theorem 34]{Bau07} are able to compute $\bbE	\EuS(X)^{(2)}_{01}$ and $\bbE\EuS(X)^{(4)}_{01}$ for $1/4<H$-fBm. Their formulae are specific to the cases $n = 2,4$ and $X$ a fBm, and are quite different to those given by \autoref{thm:main}. That the two coincide is immediate at level $2$ by \eqref{eq:ES2}, and in \autoref{appendix} we perform this check at level $4$.
	\end{expl}
	
	The following example shows how \autoref{thm:main} has the potential to generate insight into numerics of numerical schemes for rough differential equations driven by Gaussian signals.
	\begin{expl}[It\^o-Taylor expansions for solutions to RDEs driven by Gaussian signals]\label{expl:taylor0}
		Assume 
		\[
		\dif Y = V(Y) \dif \bfX, \qquad Y_0 = y_0
		\]
		is an RDE (rough differential equation) driven by the Gaussian rough path $\bfX$ (defined by the first $1$, $2$ or $3$ levels of $\EuS(X)$, depending on how rough $X$ is). Proceeding formally, and denoting by $V_{\gamma_1} \cdots V_{\gamma_n}$ composition of vector fields (and using Einstein notation), we can then expand the solution $Y$ as
		\begin{align*}
		Y_t &= \sum_{n = 0}^\infty V_{\gamma_1} \cdots V_{\gamma_n}(y_0) \EuS(X)_{0t}^{\gamma_1,\ldots,\gamma_n} \\
		&= \sum_{n = 0}^\infty V_{\gamma_1} \cdots V_{\gamma_n}(y_0) \sum_{\substack{0 \leq m \leq n \\ m \equiv n \ \mathrm{mod} \ 2}} \mathscr w^m \EuS(X)_{0t}^{\gamma_1,\ldots,\gamma_n} \\
		&= \sum_{n = 0}^\infty V_{\gamma_1} \cdots V_{\gamma_n}(y_0) \sum_{\substack{0 \leq m \leq n \\ m \equiv n \ \mathrm{mod} \ 2}} \sum_{P \in \mathcal P^n_m} \skoo^m P^{\gamma_1,\ldots,\gamma_n}_{0t} \\
		&= \sum_{m = 0}^\infty \sum_{\substack{n \geq m \\ n \equiv m \ \mathrm{mod} \ 2}}^\infty V_{\gamma_1} \cdots V_{\gamma_n}(y_0) \sum_{P \in \mathcal P^n_m} \skoo^m P^{\gamma_1,\ldots,\gamma_n}_{0t}\text{ .}
		\end{align*}
		The expansion on the first line can be viewed as the extension to the Gaussian case of Stratonovich-Taylor series, the one on the last line can be viewed as that of It\^o-Taylor series \cite{KP92}. The latter has the advantage that its terms fit in well with the Wiener chaos decomposition of $Y_t$, although it should be observed that $\mathscr w^m Y_t$ is represented as an infinite series, namely the second sum in the last line above. Also, this expansion cannot be expected to coincide with the Wiener chaos decomposition of $Y_t$ if it is performed at times other than $0$, with $Y_0 = y_0$ deterministic. This is because, unless $X$ is a martingale, the Wiener chaos isometries will not hold conditionally on $\mathcal F_s$.
	\end{expl}
	
	\begin{rem}[Stationarity and joint stationarity of increments]\label{rem:stat}
		$X$ is stationary if and only if we may write
		\begin{equation}
			R(s,t) = \overline R(t-s)
		\end{equation}
		for some function $\overline R \colon [0,T]  \to \bbR^{d\times d}$. In this case we have
		\begin{equation}
			\begin{split}
				&\partial_{12} R( s,  t) = -\overline R{}''(t-s), \quad R'( t) = 0, \quad \partial_2R(s, t) = \overline R{}'(t-s)  \\
				\implies \quad& \tfrac 12 R'(t) - \partial_2R(s,t) = - \overline R'(t-s)\text{ .}
			\end{split}
		\end{equation}
		An example of a centred stationary Gaussian process is the stationary Ornstein-Uhlenbeck process $e^{-t/2}W_{e^t}$ where $W$ is a Brownian motion and $t \in [0,T]$: its covariance function is $R(s,t) = e^{-(t-s)/2}$ for $s \leq t$. This process however, strictly speaking, is not among those considered here, as it has random initial condition.
		
		There is a much weaker property that results in a similar simplification. We will say that a stochastic process $X$ has \emph{jointly stationary increments} if for all $s_1 \leq t_1, \ldots, s_n \leq t_n $ the distribution of the random vector of increments $(X_{s_1t_1},\ldots,X_{s_nt_n})$ only depends on the differences $t_1-s_1,\ldots,t_n-s_n$ and $s_2-s_1,\ldots, s_n-s_{n-1}$ (if $n = 1$ the latter condition vanishes, and ordinary stationarity of increments is recovered). If $X$ is Gaussian this need only be required for $n = 2$, and if it holds we may write
		\begin{equation}
			R(\Delta(s,u), \Delta(t,v)) = \mathbb E[X_{su} \otimes X_{tv}] = \widehat R(u-s,v-t,t-s)
		\end{equation}
		for some function $\widehat R \colon [0,T]^3 \to \bbR^{d\times d}$. This property is satisfied by fBm, since if $H$ is the Hurst parameter we have
		\begin{align*}
			&R(\Delta(s,u), \Delta(t,v))\\
			={} &\frac 12 \big[ (t-u)^{2H} + (v-s)^{2H} - (t-s)^{2H} - (v-u)^{2H}\big] \\
			={}& \frac 12 \big[ 	\big((t-s)- (u-s)\big)^{2H} + \big((v-t) + (t-s)\big)^{2H} - \big(t-s\big)^{2H} - \big((v-t) + (t-s) - (u-s)\big)^{2H}\big]\text{ .}
		\end{align*}
		If $X$ has jointly stationary increments
		\begin{equation}
			\partial_{12}R( s,  t) = \lim_{\substack{u \to s \\ v \to t}} \frac{ R(\Delta(s,u), \Delta(t,v))}{(v-t)(u-s)} = \partial_{12}\widehat R(0,0,t-s) \text{ .}
		\end{equation}
		Although similar simplifications are not available for $\partial_2 R(s, t)$ and $R'( t)$ individually (as they are in the stationary case), they are for their difference: indeed, using that $R(\,\cdot\,,0) \equiv 0$, we have
		\begin{align*}
				& \tfrac 12 R(t+h) - \tfrac 12 R(t) - \big( R(s,t+h) - R(s,t) \big)\\
				={}& \tfrac 12 \big[ R(\Delta(s,t), \Delta(t,t+h)) + R(\Delta(s,t+h),\Delta(t,t+h)) \big]
			\end{align*}
		which implies
		\begin{align*}
			\tfrac 12 R'(t) - \partial_2R(s,t) &= \tfrac 12 \partial_h|_{h = 0} \big[\widehat R(t-s,h,t-s) + \widehat R(t+h-s,h,t-s) \big] \\
			&= \tfrac 12\partial_1\widehat R(t-s,0,t-s) + \partial_2\widehat R(t-s,0,t-s) \text{ .}
		\end{align*}
		We therefore conclude that joint stationarity of increments, though a much more general property than stationarity, results in the same simplifications that are of relevance to \autoref{thm:main}, namely that $\partial_{12}R( s,  t)$ and $\tfrac 12 R'(t) - \partial_2R(s,t)$ only depend on $t-s$. This can be of aid in simplifying the expression of the integrals in the formula for $\mathscr w^m \EuS(X)$, since it is possible to perform substitutions of the form $v_{ij} = u_j - u_i$. It does not, however, guarantee that these integrals become analytically solvable, as simple examples show (e.g.\ the integral $\int_0^1 v^{2H-1}(1-v)^{2H-1} \dif v$ appearing in \autoref{appendix}).
	\end{rem}

	We now consider a few examples of Gaussian processes to which our results apply; in all cases, $X$ will have i.i.d.\ components, and we will use $R$ to denote the scalar covariance function of each component. Arguably the most important example of a stochastic process for which the signature has not yet been computed is fractional Brownian motion in the regime of negatively-correlated increments:
	\begin{expl}[$(1/4,1/2) \ni H$-fBm]\label{expl:fBm}
		\emph{Fractional Brownian motion} with \emph{Hurst parameter} $H \in (0,1)$ ($H$\emph{-fBm}), introduced in \cite{MVN68}, is a scalar centred Gaussian process with covariance function
		\begin{equation}
			R(s,t) =  \frac {1}{2} (t^{2H} + s^{2H} - (t-s)^{2H}), \quad s \leq t\text{ .}
		\end{equation}
		It is not a semimartingale unless $H = 1/2$, in which case it is Brownian motion. Here we consider the case $H \in (1/4,1/2)$: this is well known to satisfy the preliminary hypotheses required in \autoref{sec:back}, and the smoothness conditions and bounds are simple to verify. Indeed, the integrands of interest for the formula of \autoref{thm:main} are given by ($s \leq t$)
		\begin{equation}\label{eq:t-s}
			\begin{split}
				\partial_{12}R( s,  t) &= H(2H - 1)(t-s)^{2H - 2} \\
				\tfrac 12 R'(t) - \partial_2R(s,t) &= H (t-s)^{2H - 1}\text{ .}
			\end{split}
		\end{equation}
		As predicted by \autoref{rem:stat}, these both are functions of $t-s$.
	\end{expl}

	\begin{rem}[$(1/2,1) \ni H$-fBm, \cite{Bau07}]\label{rem:fBm}
		If $R(\,\cdot\,,\,\cdot\,)$ is once differentiable on the diagonal, then
		\[
		R'(t) = \frac{\dif}{\dif t} R(t,t) = 2 \partial_2 R(t,t)
		\]  
		and we have
		\[
		\int_s^t \big[ \tfrac 12 R'(v) - \partial_2R(s,v) \big] \dif v = \int_s^t \big[ \partial_2R(v,v) - \partial_2R(s,v) \big] \dif v = \int_{s < u < v < t} \partial_{12}R(u,v) \dif u \dif v\text{ .}
		\]
		By performing this substitution in \autoref{cor:esig} for the case of $1/2 < H$-fBm (this means always applying the first case in \eqref{eq:difCases}, i.e.\ allowing arcs between consecutive nodes, which replace lines), we recover the formula of \cite[Theorem 31]{Bau07} (note that the symmetry factor --- meant to factor out permutations of pairings and transpositions within each pair --- is not present in our case, since we are summing over pairings and not permutations). Other examples of processes in a similar regularity regime are those Gaussian Volterra processes with strictly regular kernels considered in \cite{BPQ13}.
	\end{rem}
	
	The following is another example of a fractional, non-semimartingale process.
	
	\begin{expl}[The Riemann-Liouville process]\label{expl:RV}
		Another centred continuous Gaussian process, originally introduced in \cite{Levy53} and subsequently in \cite{MVN68}, is the \emph{Riemann-Liouville process} with Hurst parameter $H \in (0,1)$ (sometimes called \say{type-II fBm}), is a centred Gaussian process with covariance function \cite[p.116-117]{MarRob99}
		\begin{equation}
			\begin{split}
				R(s,t) &\underset{s<t}{=} \frac 12\bigg[t^{2H} + s^{2H} - \underbrace{2H(t-s)^{2H} \bigg( \frac{1}{2H} + \int_0^{s/(t-s)} \big((1+u)^{H-1/2} - u^{H-1/2}\big)^2\dif u  \bigg)}_{= R(\Delta(s,t),\Delta(s,t))}\bigg]\text{ .}
			\end{split}
		\end{equation}
		Like fBm, this process specifies to Brownian motion when $H = 1/2$ and is otherwise not a semimartingale. Their main difference between the two is that fBm has jointly stationary increments while for the Riemann-Liouville process not even single increments are stationary. We were not able to find a satisfactory expression for the derivatives of the covariance function of this process, and thus were not able to determine whether (for $H > 1/4$) it satisfies the conditions necessary for applying \autoref{thm:main}. However, we believe that examples such as this provide strong motivation for not confining our study to fBm and to allow for more general processes.
	\end{expl}

	Another important restriction of the main result is the following case:
	\begin{rem}[Gaussian martingales, \cite{Faw04}]\label{expl:mart}
		When $X$ is a continuous Gaussian martingale, its quadratic variation coincides with its variance function (as can be seen by the fact that $X_t^2 - R(t)$ is a martingale). The Dubins-Schwarz theorem then implies that $X$ can be represented as the deterministically-reparametrised Brownian motion $W_{R(t)}$. Assuming equal distribution of components, we can use this and the formula for the expected signature of Brownian motion \eqref{eq:sumSquares} to compute
		\begin{equation}
			\mathbb E \EuS(X)_{st}^{\gamma_1,\ldots,\gamma_{2n}} = \frac{R(\Delta(s,t))^n}{2^n n!} \kron^{\gamma_1 \gamma_2} \cdots \kron^{\gamma_{2n-1} \gamma_{2n}}\text{ .}
		\end{equation}
		Since by martingality $\partial_{12}R(s,t) = 0 = \partial_2R(s,t)$ on $s < t$, \autoref{thm:main} reduces to a sum of iterated integrals that only involve $\frac 12 R'$, which coincides with the above formula.
	\end{rem}
	
	We conclude with two examples of centred, continuous Gaussian semimartingales which are not martingales and do not have stationary increments.
	\begin{expl}[Brownian bridge returning to the origin]
		The Brownian Bridge returning to the origin at time $T$ is a process whose law is given by disintegrating the Wiener measure on the event $W_T = 0$, where $W$ is a $d$-dimensional Brownian motion starting at the origin. It can be written either as
		\[
		X_t = W_t - \frac tT W_T, \qquad t \in[0,T]
		\]
		or adaptedly as
		\[
		X_t = (T-t) \int_0^t \frac{\dif W_s}{T-s}, \qquad t \in [0,T)
		\]
		(and $X_T = 0$). Its covariance function is given by 
		\begin{equation}
			R(s,t) = s\Big( 1-\frac tT\Big), \qquad s \leq t
		\end{equation}
		and the integrands of interest are thus
		\begin{equation}\label{eq:dRBb}
			\begin{split}
				\partial_{12}R( s,  t) &= -\frac 1T \\
				\tfrac 12 R'(t) - \partial_2R(s,t) &= \frac 12 - \frac {t-s}T\text{ .}
			\end{split}
		\end{equation}
		It should be mentioned that $X$, as a process defined on $[0,T]$, fails the non-degeneracy condition \cite[p.2125]{CF10}. This is, however, not a problem, as we can view it as defined on the interval $[0,T-\varepsilon$] and obtain the signature terms $\EuS(X)_{sT}$ through a limiting argument. The bounds of \eqref{eq:requiredIneq}, which in this example and the one below only involve linear terms, are easily checked (and indeed the first is not even sharp). Note that the iterated integrals of \eqref{eq:dRBb} can all be solved explicitly as polynomials.
	\end{expl}
	
	\begin{expl}[Centred Ornstein-Uhlenbeck processes started at $0$]
		We consider an Ornstein-Uhlenbeck process with zero mean and deterministic initial condition, given by the Wiener-It\^o integral
		\[
		X_t = \sigma\int_0^t e^{-\theta(t-u)} \dif W_u
		\]
		with $\sigma,\theta \in (0,+\infty)$. Its covariance function is given by
		\begin{equation}
			R(s,t) = \frac{\sigma^2}{2\theta} \big( e^{-\theta(t-s)} - e^{-\theta(s+t)} \big), \qquad s \leq t
		\end{equation}
		and $\partial_{12}R(\dif s, \dif t)$, $\tfrac 12 R'(t) - \partial_2R(s,t)$ can be computed directly. Once again, all conditions are satisfied (see \cite[p2138]{CF10}).
	\end{expl}

	\section{Proof of the main result}\label{sec:proof}
	
	Recall that we are using $\lesssim$ to denote inequalities whose constant of proportionality may only depend on $T,H$ and other properties of a fixed process $X$. Since most of the arguments presented in this section only concern bounds and convergence, we will suppress indices (i.e.\ treat the scalar case) most of the time, so as not to clutter the notation. Given $P \in \mathcal P^n_m$, denote $|P|_{st}$ the function $\Delta^m[s,t] \to \bbR$ defined analogously to \autoref{def:defIntegral}, but replacing each integrand $\partial_{12}R(u,v)$ with $(v-u)^{2H-2}$ and each integrand $\frac 12 R'(v) - \partial_2 R(u,v)$ with $(v-u)^{2H-1}$. For example, if $Q$ is the diagram of \eqref{eq:diagEx}
	\[
	|Q|_{st} = \mathbbm 1_{\Delta^2[s,t]}(u_2,u_7)  \int_{\substack{s < u_1 < u_2 \\ u_2 < u_3 < u_4 < u_6 < u_7 \\ u_7 < u_8 < t}} (u_4 - u_1)^{2H-2} (u_8 - u_3)^{2H-2} (u_6 - u_4)^{2H-1}\dif u_1 \dif u_3 \dif u_4 \dif u_6 \dif u_8\text{ .}
	\]
	The following proposition guarantees that all the integrals considered in the main theorem are convergent.
	\begin{prop}[Finite improper integrals]\label{prop:finite}
		For $m \leq n$ and $P \in \mathcal P^n_m$
		\begin{equation}
			|P|_{st}  \lesssim (t-s)^{(n-m)H}
		\end{equation}
		uniformly over $\Delta^m[s,t]$.
		\begin{proof}
			We proceed by induction on $n-m$. When $P$ only has single nodes ($m = n$) the statement is trivial. We will proceed by considering several cases for the last node in $P$; the simplest of these occurs when it is single: the statement follows immediately from the inductive hypothesis. For the next case, we will need the following bound:
			\begin{equation*}
				\begin{split}
					&\int_{\Delta^n[s,t]} (u_1 - s)^{2H - 1} \cdots (u_n - u_{n-1})^{2H - 1} \dif u_1 \cdots \dif u_n \\
					\lesssim{}&\int_{\Delta^{n-1}[s,t]} (u_1 - s)^{2H - 1} \cdots (u_{n-1} - u_{n-2})^{2H - 1} (t - u_{n-1})^{2H} \dif u_1 \cdots \dif u_{n-1}\\
					\leq{}&(t-s)^{2H}\int_{\Delta^{n-1}[s,t]} (u_1 - s)^{2H - 1} \cdots (u_{n-1} - u_{n-2})^{2H - 1} \dif u_1 \cdots \dif u_{n-1} \\
					\lesssim{}& \ldots \lesssim (t-s)^{2nH}\text{ .}
				\end{split}
			\end{equation*}
			For a diagram $C$ whose last node is the right endpoint of an arc, using the bound above we have
			\begin{equation*}
				\begin{split}
					&|C\,\tikz[baseline = -0.1ex]{
						\draw [decorate,decoration={brace,amplitude=4pt,mirror,raise=0.5ex}]
						(-0.1 ,0) -- (1,0);
						\draw[fill] (0,0) circle [radius=0.04];
						\draw[fill] (0.9,0) circle [radius=0.04];
						\draw (0,0)-- (0,0.3);
						\draw (0.9,0)-- (0.9,0.3);
						\node at (0.45,0){$\ldots$};
						\node at (0.45,-0.35){$\scriptstyle n$};
					}|_{st}\\
					={}&\int_{\Delta^{n+1}[s,t]}  |C|'_{su_0}        
					(u_1-u_0)^{2H-1} \cdots (u_n-u_{n-1})^{2H-1} \dif u_0 \cdots \dif u_n \\
					= {}& \int_s^t |C|'_{su_0}  \int_{\Delta^n[u_0,t]}         
					(u_1-u_0)^{2H-1} \cdots (u_n-u_{n-1})^{2H-1} \dif u_1 \cdots \dif u_n \ \, \dif u_0 \\
					\lesssim{} &\int_s^t |C|'_{su_0} (t-u_0)^{2nH} \dif u_0 \\
					\leq{} &(t-s)^{2nH}\int_s^t |C|'_{su_0} \dif u_0 \\
					\leq{} &(t-s)^{2nH} |C|_{st}
				\end{split}
			\end{equation*}
			where $|C|_{su_0}'$ equals the integral representing $|C|_{su_0}$ with the only difference that we are not integrating w.r.t.\ the variable $u_0$ in $(u_0-r)^{2H-2}$, which represents the arc that terminates at the last node of $C$. Similarly, if the last node in $C$ is single, we have
			\begin{equation*}
				\begin{split}
					&|C\,\tikz[baseline = -0.1ex]{
						\draw [decorate,decoration={brace,amplitude=4pt,mirror,raise=0.5ex}]
						(-0.1 ,0) -- (1,0);
						\draw[fill] (0,0) circle [radius=0.04];
						\draw[fill] (0.9,0) circle [radius=0.04];
						\draw (0,0)-- (0,0.3);
						\draw (0.9,0)-- (0.9,0.3);
						\node at (0.45,0){$\ldots$};
						\node at (0.45,-0.35){$\scriptstyle n$};
					}|_{st}\\
					={}&|C|_{su_0} \int_{\Delta^n[u_0,t]}       
					(u_1-u_0)^{2H-1} \cdots (u_n-u_{n-1})^{2H-1} \dif u_1 \cdots \dif u_n\\ \lesssim{} &|C|_{su_0} (t-s)^{2nH} \\ \leq{} &|C|_{st} (t-s)^{2nH}
				\end{split}
			\end{equation*}
			where $C$ is not differentiated since it terminates in a node representing a free variable, $u_0$. We now consider arcs: assume there are $i$ arcs/lines within $A$, $j$ within $B$, and that there are $k$ arcs between nodes in $A$ and nodes in $B$ (collectively represented below by the dashed arc). Let $A^\circ$ and $B^\circ$ denote the diagrams given by eliminating such arcs from $A$ and $B$: the nodes that have become single as a result now represent free variables, which we call $w_1,\ldots,w_k$, $z_1,\ldots,z_k$. We first consider the case in which $j>0$:
			\begin{equation*}
				\begin{split}
					&|\!\!\tikz[baseline = -0.7ex]{
						\node at (0,0) {$A$};
						\draw[fill] (0.4,-0.1) circle [radius=0.035];
						\node at (0.8,0) {$B$};
						\draw[fill] (1.2,-0.1) circle [radius=0.035];
						\draw (0.4,-0.1) .. controls (0.4,0.5) and (1.2,0.5) .. (1.2,-0.1);
						\draw[dashed] (0,0.2) .. controls (0,0.5) and (0.8,0.5) .. (0.8,0.2);
					} \, |_{st} \\
					={}& \int_{s < u < v < t} \int_{\Delta^k[s,u] \times \Delta^k[u,v]} |A^\circ|_{su} | B^\circ|_{uv}  (z_1 - w_1)^{2H-2} \cdots (z_k - w_k)^{2H-2} \dif w_1 \dif z_1 \cdots \dif w_k \dif z_k \\ &\quad(v-u)^{2H-2} \dif u \dif v \\
					\lesssim{}& \int_{s < u < v < t} \int_{\Delta^k[s,u] \times \Delta^k[u,v]} (u-s)^{2iH} (v-u)^{2jH}  (z_1 - w_1)^{2H-2} \cdots (z_k - w_k)^{2H-2} \dif w_1 \dif z_1 \cdots \dif w_k \dif z_k \\ &\quad(v-u)^{2H-2} \dif u \dif v \\
					\leq{} &(t-s)^{2iH} \int_{s < u < v < t} (v-u)^{2jH + 2H - 2}  \int_{[s,u]^k \times [u,v]^k} ( z_1 - w_1)^{2H-2} \cdots ( z_k - w_k)^{2H-2}  \dif w_1 \dif z_1\cdots \\ &\quad \cdots\dif w_k \dif z_k \ \dif u \dif v \\
					\lesssim{}& (t-s)^{2iH} \int_{s < u < v < t} (v-u)^{2(j+1)H-2} | (v-u)^{2H} + (u-s)^{2H} - (v-s)^{2H} |^k \dif u \dif v \\
					\lesssim{}& (t-s)^{2(i+k)H} \int_{s < u < v < t} (v-u)^{2(j+1)H-2} \dif u \dif v \\
					\leq{}&(t-s)^{2(i+j+k)H}
				\end{split}
			\end{equation*}
			where we have used $2H(j+1)-1 \geq 4H-1 >0$ since $H >1/4$. Note that the absolute values in the third-last expression can be removed by separately considering the cases $H >1/2$ and $H < 1/2$. Assume instead $j=0$: this means $B$ must contain at least one node that is either single or paired with a node in $A$; it cannot be that $B = \varnothing$ or the diagram would contain an arc between two consecutive nodes, which is ruled out. The case in which there is a node in $B$ which is single (see \autoref{fig:gradient}) does not require $H > 1/4$: letting $r$ denote the free variable represented by such a node, and proceeding similarly to the above, we have
			\begin{equation*}
				\begin{split}
					&|\!\!\tikz[baseline = -0.7ex]{
						\node at (0,0) {$A$};
						\draw[fill] (0.4,-0.1) circle [radius=0.035];
						\node at (0.8,0) {$B$};
						\draw[fill] (1.2,-0.1) circle [radius=0.035];
						\draw (0.4,-0.1) .. controls (0.4,0.5) and (1.2,0.5) .. (1.2,-0.1);
						\draw[dashed] (0,0.2) .. controls (0,0.5) and (0.8,0.5) .. (0.8,0.2);
					}\,|_{st}\\
				={}&\int_{s < u < r < v < t} \int_{\Delta^k[s,u] \times \Delta^k[u,v]} |A^\circ|_{su} (z_1 - w_1)^{2H-2} \cdots (z_k - w_k)^{2H-2} \dif w_1 \dif z_1 \cdots \dif w_k \dif z_k \\ &\quad(v-u)^{2H-2} \dif u \dif v \\
					\lesssim{}& (t-s)^{2(i+k)H} \int_{s < u < r < v < t} (v-u)^{2H-2} \dif u \dif v \\
					\lesssim{} &(t-s)^{2(i+k)H} | (t-r)^{2H} + (r-s)^{2H} - (t-s)^{2H} | \\
					\lesssim{}& (t-s)^{2(i+k+1)H}\text{ .}
				\end{split}
			\end{equation*}
			Finally, consider the case in which $j = 0$ and $k > 0$ (and $B$ may have no single nodes):
			\begin{equation*}
				\begin{split}
					&|\!\tikz[baseline = -0.7ex]{
						\node at (0,0) {$A$};
						\draw[fill] (0.4,-0.1) circle [radius=0.035];
						\node at (0.8,0) {$B$};
						\draw[fill] (1.2,-0.1) circle [radius=0.035];
						\draw (0.4,-0.1) .. controls (0.4,0.5) and (1.2,0.5) .. (1.2,-0.1);
						\draw[dashed] (0,0.2) .. controls (0,0.5) and (0.8,0.5) .. (0.8,0.2);
					}\, |_{st}\\
					\lesssim{} &(t-s)^{2iH} \int_{s  < u  < v < t} \int_{[s,u]^k \times \Delta^k[u,v]} ( z_1 - w_1)^{2H-2} \cdots ( z_k - w_k)^{2H-2}   \dif w_1\dif z_1 \cdots\dif w_k \dif z_k \\ &\quad(v-u)^{2H-2} \dif u  \dif v \\
					\lesssim{}& (t-s)^{2iH}\int_{s < u < z_1 < \ldots < z_k < t}   |(z_1 - u)^{2H-1} - (z_1 - s)^{2H-1}| \cdots |(z_k - u)^{2H-1} - (z_k - s)^{2H-1}| \\ &\quad|(t-u)^{2H-1} - (z_k-u)^{2H-1}| \dif u \dif z_1 \cdots \dif z_k \\
					\lesssim{}& (t-s)^{2(i+k-1)H}\int_{s < u < z_k < t} |(z_k - u)^{2H-1} - (z_k - s)^{2H-1}| |(t-u)^{2H-1} - (z_k-u)^{2H-1}| \dif u \dif z_k\text{ .}
				\end{split}
			\end{equation*}	
			Once again, the absolute values distinguish between $H \lessgtr 1/2$. Expanding the product, we observe that three of the integrals feature products of different terms, each to the power of $2H-1$: in these, at least one of $z_k$ or $u$ only appears once, which means this variable may be integrated out and the resulting term bounded (up to a constant) by $(t-s)^{2H}$, with the remaining integral solved similarly. The fourth integral instead is $\int_{s < u < z < t} (z-u)^{4H-2} \dif u \dif z$ which is finite again thanks to $H > 1/4$. This shows that we have $\lesssim (t-s)^{2(i+k+1)H}$ in the above expression and concludes the proof.
		\end{proof}
	\end{prop}
\begin{rem}[Modified $|P|$]\label{rem:finite}
We have stated the previous proposition under in the most natural manner; in particular note how, in the prototypical case of fBm, the integrals $|P|_{st}$ are multiples of $P_{st}$. We will, however, additionally need a slightly modified version of this result, in which the definition of $|P|$ is changed as follows: maximal sequences 
\[
\int_{\Delta^k[u,v]} (w_1 - u)^{2H - 1} \cdots (w_k - w_{k-1})^{2H - 1} \dif w_1 \cdots \dif w_k
\]
occurring in the middle of the expression for $|P|$, are replaced with their bound $(v-u)^{2kH}$, and each integrand $(v-u)^{2H-2}$ is replaced with $((v-u) \wedge 1/2)^{2H-2}$. That the statement continues despite these modifications to hold is obvious for the first, and for the second it follows from the facts that all integrals are still convergent (by the same proof) and the $1/2$ can be replaced with $1/2 \wedge T$ and absorbed in the constant of proportionality.
\end{rem}
	
Just like in \cite{Bau07}, we approximate $X$ piecewise linearly. Let $X^\ell$ be a sequence of piecewise linear approximations of $X$ along partitions $\pi_\ell$ on $[0,T]$ with step size that vanishes as $\ell \to \infty$. It will be helpful to assume that the intervals in the mesh $\pi_\ell$ all have the same length $\varrho_\ell$; this simplifying assumption can be made because it is only necessary to show convergence along \emph{a} sequence of such approximations, since it is known that the limit does not depend on the particular choice of $\pi_\ell$ (or indeed on the type of piecewise smooth approximation in a broad class of these) \cite[Ch.\ 15]{FV10}. For $t \in[0,T]$ we will write $t_\ell^-$ and $t_\ell^+$ to respectively denote the endpoints $a$ and $b$ of the interval of $\pi_\ell$ s.t.\ $t \in [a,b)$. Explicitly, $X^\ell$ and its piecewise-defined derivative are given by
\begin{equation}\label{eq:pwl}
	\begin{split}
	X^\ell_t &= X_{t^-_\ell} + \varrho^{-1}_\ell (t-t^-_\ell) X_{t^-_\ell t^+_\ell} \\
	\dot X^\ell_t &= \varrho^{-1}_\ell X_{t^-_\ell t^+_\ell}
	\end{split}
\end{equation}
where, as usual, $X_{ab} \coloneqq X_b - X_a$ denotes the increment. In order to use Stroock's formula \eqref{eq:stroock}, we will be considering Malliavin derivatives of the signature of the piecewise-linear interpolations of $X$, 
\[
\EuS(X^\ell)^{\gamma_1,\ldots,\gamma_n}_{st} = \int_{\Delta^n[s,t]} \dot X^{\ell;\gamma_1}_{u_1}  \cdots  \dot X^{\ell;\gamma_n}_{u_n} \dif u_1 \cdots \dif u_n \text{ ,}
\]
which in turn requires us to consider those of the single factors:
\begin{equation}\label{eq:Dpwl}
	\EuD_{v} \dot X_u^{\ell;\gamma} = \varrho^{-1}_\ell  \mathbbm 1^\gamma_{[u^-_\ell,u^+_\ell)}(v) = \varrho^{-1}_\ell  \mathbbm 1^\gamma_{[v^-_\ell,v^+_\ell)}(u)\text{ .}
\end{equation}
For $P \in \mathcal P^n_m$, we provide a discretised analogue to \autoref{def:defIntegral}:
\begin{defn}[$P^{\ell;\gamma_1,\ldots,\gamma_n}_{st}$]\label{defn:seql}
	For $\gamma_1,\ldots,\gamma_n \in [d]$, $0 \leq s \leq t \leq T$ and $P \in \mathcal P^n_m$ define
	\begin{equation}
		\begin{split}
			P_{st}^{\ell;\gamma_1,\ldots,\gamma_n}(v_k \mid k \in \overline P)
			&\coloneqq  \int_{\Delta^{n}[s,t]} \prod_{\{i,j\} \in P} \bbE [\dot X^{\ell;\gamma_i}_{u_i} \dot X^{\ell;\gamma_j}_{u_j}] \dif u_i \dif u_j  \cdot  \prod_{k \in \overline P} \varrho_\ell^{-1} \mathbbm 1^{\gamma_k}_{[v^-_{k;\ell}v^+_{k;\ell})}(u_k)\dif u_k\text{ .}
		\end{split}
	\end{equation}
as an element of $\mathcal E^{\otimes m}$, whose arguments are given to the functions $\mathbbm 1^{\gamma_k}_{[u^-_{k;\ell}u^+_{k;\ell})}$ with $k \in \overline P$.
\end{defn}
Note how the above definition, unlike \autoref{def:defIntegral} does not distinguish between consecutive and non-consecutive pairings: this will only become important in the limit. Moreover, we are integrating over all $n$ variables, including the $u_k$ with $k \in \overline P$: this is because the time arguments of the function, $v_k$, are supplied separately, with the respective index variables supplied as arguments to $\kron_{\gamma_k}$, $k \in \overline P$. The functions $P^\ell_{st}$ are summands in the expression of which we want to compute the limit:
\begin{lem}[Expected Malliavin derivatives of signature approximations]\label{lem:expMal}
\begin{equation}
	\bbE \EuD^m \EuS(X^\ell)_{st}^{\gamma_1,\ldots,\gamma_n} = m!\sum_{P \in \mathcal P^n_m} P_{st}^{\ell;\gamma_1,\ldots,\gamma_n} \ \in \mathcal E^{\odot m}
\end{equation}
\begin{proof}
	This is a consequence of \eqref{eq:pwl}, \eqref{eq:Dpwl}, the (iterated) Leibniz rule for the Malliavin derivative and Wick's formula for the mixed moments of a Gaussian vector (as it was already used in \cite[Theorem 31]{Bau04}). The details are a matter of simple combinatorics; in particular note how, instead of summing over $m!$ terms corresponding to the ways of permuting the $m$ derivatives (for a fixed $P \in \mathcal P^n_m$), we are only including the term corresponding to the identity permutation and multiplying by $m!$, which identifies the same element of $\mathcal E^{\otimes m}$ up to symmetry.
\end{proof}
\end{lem}

In order to prove convergence, it is unfortunately not possible to argue by dominated convergence applied to \autoref{defn:seql}: this is because the factors in the integrand given by consecutive pairings $\bbE [\dot X^{\ell;\gamma_i}_{u_i} \dot X^{\ell;\gamma_{i+1}}_{u_{i+1}}]$ converge to non-integrable functions (e.g.\ $(v-u)^{2H-2}$ on $\Delta^2[s,t]$ for fBm) and the ones corresponding to Malliavin derivatives $\varrho_\ell^{-1} \mathbbm 1^{\gamma_k}_{[v^-_{k;\ell}v^+_{k;\ell})}(u_k)\dif u_k$ do not converge at all (in fact they converge, as distributions, to Dirac deltas $\kron_{v_k}$). The reason that convergence holds is that all these quantities are integrated. To successfully exploit this, we will write each integral $P^\ell_{st}$ as a nested integral, distinguishing between the three types of integrands:
	\begin{equation}\label{eq:3ints}
	\int (\text{non-consecutive pairings})  \int (\text{Malliavin derivatives}) \prod_{\substack{\text{maximal} \\ \text{sequences}}} \int (\text{consecutive pairings})\text{ .}
\end{equation}
The outer integral contains the product of all terms $\bbE [\dot X^{\ell;\gamma_i}_{u_i} \dot X^{\ell;\gamma_j}_{u_j}]$ with $|j-i| > 1$. These are multiplied with the second integral, which integrates all factors coming from Malliavin derivatives. Finally, we partition the remaining integrands $\bbE [\dot X^{\ell;\gamma_h}_{u_h} \dot X^{\ell;\gamma_{h+1}}_{u_{h+1}}]$ into maximal sequences and integrate each individually: these integrals are integrands in the second integral, alongside the Malliavin derivatives. The operations of exchanging the order of integrals are all justified by Fubini's theorem, considering that all integrals are actually finite sums. We illustrate all of this with a simple example: consider the diagram (suppressing indices)
\[
P \coloneqq \tikz[baseline = -0.1ex]{
	\draw[fill] (0,0) circle [radius=0.04];
	\draw[fill] (0.3,0) circle [radius=0.04];
	\draw[fill] (0.6,0) circle [radius=0.04];
	\draw[fill] (0.9,0) circle [radius=0.04];
	\draw[fill] (1.2,0) circle [radius=0.04];
	\draw (0,0) .. controls (0,0.4) and (1.2,0.4) .. (1.2,0);
	\draw (0.9,0)-- (0.9,0.2);} \ \in \mathcal P^6_2\text{ .}
\]
According to \autoref{defn:seql}, we have 
\[
P^\ell_{st}(v_1,v_2) = \int_{\Delta^6[s,t]} \varrho_\ell^{-2} \mathbbm 1_{[v^-_{2;\ell},v^+_{2;\ell})}(u_2) \mathbbm 1_{[v^-_{3;\ell},v^+_{3;\ell})}(u_3) \bbE[\dot X^\ell_{u_1} \dot X^\ell_{u_6}] \bbE[\dot X^\ell_{u_4} \dot X^\ell_{u_5}] \dif u_1 \dif u_4 \dif u_5 \dif u_6\text{ .}
\]
Re-organising this expression as described in \eqref{eq:3ints} we obtain
\begin{align*}
&\int_{s < u_1 < u_6 < t} \bbE[\dot X^\ell_{u_1} \dot X^\ell_{u_6}]  \\
&\qquad\bigg[\int_{u_1 < u_2 < u_3 <u_6} \varrho_\ell^{-2} \mathbbm 1_{[v^-_{2;\ell},v^+_{2;\ell})}(u_2) \mathbbm 1_{[v^-_{3;\ell},v^+_{3;\ell})}(u_3) \\
&\qquad\qquad\bigg[\int_{u_3 < u_4 < u_5 < u_6} \bbE[\dot X^\ell_{u_4} \dot X^\ell_{u_5}]\dif u_4 \dif u_5 \bigg] \dif u_2 \dif u_3 \bigg] \dif u_1 \dif u_6\text{ .}
\end{align*}
Note that the domain of integration of the innermost integral can be described in terms of variables of the two outer integrals: this extends to the case in which there is more than one maximal sequence, by maximality, and is crucial for the factorisation into integrals over maximal sequences to be possible.

The reason for the nested rewriting of \eqref{eq:3ints} is that it will be possible to show convergence of the integrals over maximal sequences, then by a separate argument infer the convergence of the middle integral, and finally by dominated convergence conclude that the outer integrals converge. We preface the proof of convergence with a few lemmas; the first of these considers the case of a single consecutive pairing, and will form the base case of an induction that handles maximal sequences of arbitrary length.

	\begin{lem}[One consecutive pairing]\label{lem:onePair}
	\[\lim_{\ell \to \infty} \int_{s < u < v < t} \mathbb E[\dot X_u^\ell \dot X_v^\ell] \dif u\dif v = \frac 12 \mathbb E[X^2_{st}] = \int_s^t \big[\tfrac 12 R'(v) - \partial_2R(s,v) \big] \dif v \]
	and the convergents are uniformly bounded by $\lesssim (t-s)^{2H}$.
	\begin{proof}
		Considering that $\dot X^\ell$ is a piecewise-constant, and that the integral on the right is therefore a finite sum, we can write
		\begin{align*}
			\int_{s < u < v < t} \bbE[\dot X^\ell_u  \dot X^\ell_v] \dif u \dif v
			={}	&\mathbb E \int_{s < u < v < t} \dot X^\ell_u  \dot X^\ell_v \dif u \dif v\\ ={}&\frac 12 \mathbb E[(X_{st}^\ell)^2] \\ \xrightarrow{\ell \to \infty} {}&\frac 12 \mathbb E[X_{st}^{ 2}] \\
			={} &\frac 12 R(\Delta(s,t), \Delta(s,t)) \\
			={}&\frac{R(s) + R(t)}{2} - R(s,t) \\
			={} &\int_s^t \big[\tfrac 12 R'(v) - \partial_2R(s,v) \big] \dif v
		\end{align*}
		where we have used that $(X^\ell_{st})^2 \xrightarrow{\ell \to \infty} X_{st}^2$ in $L^2$. For the second statement, we rely on the first two identities above and distinguish between the cases $s^-_\ell = t^-_\ell$ and $s^-_\ell < t^-_\ell$: in the former we have, using \eqref{eq:onDiagEst}
		\begin{align*}
			|\mathbb E[(X^\ell_{st})^2]| = |\mathbb E[\varrho_\ell^{-2} (t-s)^2 X_{s^-_\ell s^+_\ell}^2]| \lesssim \varrho_\ell^{2H-2} (t-s)^2 \leq (t-s)^{2H}
		\end{align*}
		since
		\[
		\Big(\frac{t-s}{\varrho_\ell} \Big)^{2-2H} \leq 1
		\]
		by $H < 1$ and $t-s \leq \varrho_\ell$. Let now $s^-_\ell < t^-_\ell$:
		\begin{align*}
			|\mathbb E[(X^\ell_{st})^2]| &= |\mathbb E[(X^\ell_{ss^+_\ell} + X^\ell_{s^+_\ell t^-_\ell} + X^\ell_{t^-_\ell t})^2]| \\
			&= \Big|\mathbb E\Big[ \big(\varrho^{-1}(s^+_\ell - s)X^\ell_{s s^+_\ell} + X_{s^+_\ell t^-_\ell} + \varrho^{-1}(t-t^-_\ell) X_{t^-_\ell t}\big)^2 \Big]\Big| \\
			&\lesssim \varrho^{-2}(s^+_\ell - s)^2 \mathbb E[X_{ss^+_\ell}^2] + \mathbb E[X_{s^+_\ell t^-_\ell}^2] + \varrho^{-2} (t-t^-_\ell)^2 \bbE[X_{t^-_\ell t}^2] \\
			&\lesssim (s^+_\ell - s)^{2H} + (t^-_\ell - s^+_\ell)^{2H} + (t-t^-_\ell)^{2H} \\
			&\lesssim (t-s)^{2H}
		\end{align*}
		by $\mathscr l^2$-Jensen's inequality, the previous case, and again \eqref{eq:onDiagEst}.
	\end{proof}
\end{lem}

The case of several consecutive pairings is more difficult to handle, and in \autoref{prop:conv} convergence of these terms will be bootstrapped from terms that only contain shorter sequences of consecutive pairings, and the above single case, by means of an inductive argument. It is worth remarking that the plausible strategy of handling these integrands together with the others by integrating only one of the variables fails:
\begin{rem}[Lack of convergence of $\mathbb E[X_{uv}^\ell \dot X_v^\ell{]}$]
	One way of dealing with sequences of consecutive pairings is by rewriting them as
	\begin{equation}\label{eq:consec}
		\begin{split}
			&\int_{s < u_1 < v_1 < \ldots < u_n < v_n < t} \mathbb E[\dot X^\ell_{u_1} \dot X^\ell_{v_1}]\cdots \mathbb E[\dot X^\ell_{u_n} \dot X^\ell_{v_n}] \dif u_1 \dif v_1 \cdots \dif u_n \dif v_n \\ ={} &\int_{\Delta^n[s,t]} \mathbb E[X^\ell_{sv_1} \dot X^\ell_{v_1}] \mathbb E[X^\ell_{v_1 v_2} \dot X^\ell_{v_2}] \cdots \mathbb E[X^\ell_{v_{n-1} v_k} \dot X^\ell_{v_n}] \dif v_1 \cdots \dif u_n\text{ .}
		\end{split}
	\end{equation}
	This has the benefit of expressing the convergents as integrals over $n$, and not $2n$, variables. The problem with this strategy is that it does not hold that $\mathbb E[X_{uv}^\ell \dot X_v^\ell] \xrightarrow{\ell \to \infty} \frac 12 R'(v) - \partial_2R(u,v)$: a simple calculation reveals
	\begin{align*}
		&\mathbb E[X_{uv}^\ell \dot X_v^\ell] \\
		={}& \varrho^{-1}_\ell \Big[ \big(1-\varrho^{-1}_\ell(v-v^-_\ell)\big)R(v^-_\ell, \Delta(v^-_\ell, v^+_\ell)) + \varrho^{-1}_\ell(v-v^-_\ell) R(v^+_\ell,\Delta(v^-_\ell,v^+_\ell))  \Big] \\
		&-\varrho^{-1}_\ell \Big[ \big(1-\varrho^{-1}_\ell(u-u^-_\ell)\big)R(u^-_\ell, \Delta(v^-_\ell, v^+_\ell)) + \varrho^{-1}_\ell(u-u^-_\ell) R(u^+_\ell,\Delta(v^-_\ell,v^+_\ell))  \Big]\text{ .}
	\end{align*}
	While the second term converges to $\partial_2R(u,v)$ (e.g.\ by the intermediate value theorem applied on the interval $[v^-_\ell,v^+_\ell]$), the first does not converge in general. To see why, it suffices to take $X$ to be Brownian motion and $\pi_\ell$ to by a diadic sequence: the first term on the right above is then equal to $\varrho^{-1}_\ell(v-v^-_\ell)$ which is indeterminate in view of the fact that for $v$ in a set of full Lebesgue measure its decimal expansion contains infinitely many $00$'s and $11$'s. The fractional case with $H < 1/2$ appears even worse behaved, i.e.\ divergent in a possibly indeterminate fashion.
\end{rem}

We now move outward in \eqref{eq:3ints} and prove a lemma that will guarantee convergence of the middle integral, conditional on the convergence of the inner ones.

\begin{lem}\label{lem:middle}
Let $f_\ell \colon [0,T]^m \to \bbR$ be a uniformly bounded sequence of functions that are continuous and piecewise smooth on the mesh $\pi_\ell$. Assume that $f_\ell$ converges to $f \colon [0,T]^m \to \bbR$ uniformly. Then
\begin{equation*}
	\begin{split}
\int_{\Delta^m[s,t]} f_\ell(u_1,\ldots,u_m) \varrho^{-m}_\ell \prod_{k = 1}^m \mathbbm 1_{[v^-_{k;\ell},v^+_{k;\ell})}(u_k) \dif u_1,\ldots,\dif u_m
 \xrightarrow{\ell \to \infty}\mathbbm 1_{\Delta^m[s,t]}(v_1,\ldots,v_m) f(v_1,\ldots,v_m)
\end{split}
\end{equation*}
where the convergence is a.e.\ in the variables $(v_1,\ldots,v_m) \in [0,T]^m$. Moreover, the convergents are uniformly bounded by $\sup_\ell\lVert f_\ell \rVert_\infty$.
\begin{proof}
	The second statement holds by uniform boundedness of $f_\ell$ and the fact that
	\[
	\int_{\Delta^m[s,t]} \varrho^{-m}_\ell \prod_{k = 1}^m \mathbbm 1_{[v^-_{k;\ell},v^+_{k;\ell})}(u_k) \dif u_1,\ldots,\dif u_n \leq 1 \text{ .}
	\]
	We will prove pointwise convergence on the subset
	\[
	[0,T]^m_* \coloneqq \{(v_1,\ldots,v_m) \in [0,T]^m \mid v_i \neq v_j \text{ for } i \neq j\}
	\]
	of $[0,T]^m$ of full Lebesgue measure. For $(v_1,\ldots,v_m) \in [0,T]^m_*$, we may, without loss of generality, start the sequence when $\ell$ is already large enough so that $[v^-_{i;\ell},v^+_{i;\ell}) \cap [v^-_{j;\ell},v^+_{j;\ell}) = \varnothing$ for $i \neq j$, where we are including $v_0 \coloneqq s$ and $v_{m+1} \coloneqq t$ in this requirement. By the mean value theorem applied individually to each $u_k$, there exist $w_{k;\ell} \in (v^-_{k;\ell},v^+_{k;\ell})$ s.t.\ 
	\[
\int_{\Delta^m[s,t]} f_\ell(u_1,\ldots,u_m) \varrho^{-m}_\ell \prod_{k = 1}^m \mathbbm 1_{[v^-_{k;\ell},v^+_{k;\ell})}(u_k) \dif u_1,\ldots,\dif u_m = \mathbbm 1_{\Delta^m[s,t]}(v_1,\ldots,v_m) f_\ell(w_{1;\ell},\ldots,w_{m;\ell}) 
	\]
	 and
	\begin{align*}
		&|\mathbbm 1_{\Delta^m[s,t]}(v_1,\ldots,v_m) f_\ell(w_{1;\ell},\ldots,w_{m;\ell}) - \mathbbm 1_{\Delta^m[s,t]}(v_1,\ldots,v_m)f(v_1,\ldots,v_m)|\\
		\leq{}&\mathbbm 1_{\Delta^m[s,t]}(v_1,\ldots,v_m) \big[|f_\ell(w_{1;\ell},\ldots,w_{m;\ell}) - f(w_{1;\ell},\ldots,w_{m;\ell})| \\
		&+|f(w_{1;\ell},\ldots,w_{m;\ell}) - f(v_1,\ldots,v_m)|\big] \\
		\leq{} &\mathbbm 1_{\Delta^m[s,t]}(v_1,\ldots,v_m) \big[\lVert f_\ell - f \rVert_\infty  + \omega^f_{(v_1,\ldots,v_m)}(\varrho_\ell) \big] 
	\end{align*}
	where $\omega^f_{(v_1,\ldots,v_m)}$ is the modulus of continuity of $f$ at the point $(v_1,\ldots,v_m)$. Both summands on the right hand side above vanish in the limit of $\ell \to \infty$, the first by uniform convergence and the second by continuity of the uniform limit of continuous functions.
\end{proof}
\end{lem}

The next two results constitute the core of our argument. They both rely on the same induction used to reduce the length of consecutive pairings, the base case of which is provided by \autoref{lem:onePair}. To illustrate it at level $4$, letting $Y$ be a stochastic process (which below will be taken to be $X^\ell$ and $X$) we have for $\alpha \neq \beta$
\begin{align*}
2 \mathbb E \EuScript{S}(Y)^{\alpha\alpha\beta\beta}_{st} = \mathbb E\EuScript{S}(Y)^{\alpha\alpha}_{st} \cdot \mathbb E \EuScript{S}(Y)^{\beta\beta}_{st} - 2 \mathbb E \EuScript{S}(Y)^{\alpha\beta\alpha\beta}_{st} - 2\mathbb E \EuScript{S}(Y)^{\alpha\beta\beta\alpha}_{st}
\end{align*}
by the shuffle property \eqref{eq:shuffle}, using identical distribution of components to group together $2\mathbb E \EuScript{S}(Y)^{\alpha\alpha\beta\beta}_{st} = \mathbb E \EuScript{S}(Y)^{\alpha\alpha\beta\beta}_{st} + \mathbb E \EuScript{S}(Y)^{\beta\beta\alpha\alpha}_{st}$ (and similar on the right hand side), and using independence of components to write $\mathbb E[\EuScript{S}(Y)_{st}^{\alpha\alpha}\EuScript{S}(Y)_{st}^{\beta\beta}] = \mathbb E\EuScript{S}(Y)^{\alpha\alpha}_{st} \cdot \mathbb E \EuScript{S}(Y)^{\beta\beta}_{st}$. While the left hand side contains a sequence of two consecutive pairs, only sequences of consecutive pairs of length one appear on the right.

\begin{lem}[Dominating function]\label{lem:dom}
For $P \in \mathcal P^n_m$ it holds that the integrand of the outermost integral of $P^\ell_{st}$ expressed in the nested form \eqref{eq:3ints}, is absolutely bounded by an integrable function, uniformly in $\ell$ and on $\Delta^m[s,t]$, so that $|P^\ell_{st}|\lesssim (t-s)^{2(n-m)H'}$ for any $1/4 < H' < H$. 
\begin{proof}
	
We begin by bounding expectations corresponding to non-consecutive pairings. As done in the proof of \cite[Theorem 31]{Bau07}, we now consider the terms
	\[
	\mathbb E[\dot X_u^\ell \dot X^\ell_v] = \varrho^{-2}_\ell R(\Delta(u^-_\ell, u^+_\ell),\Delta(v^-_\ell, v^+_\ell))
	\]
	in three different cases: for $u^-_\ell = v^-_\ell$
	\[
	|\mathbb E[\dot X^\ell_u \dot X^\ell_v]|  \lesssim \varrho_\ell^{2H-2} \leq (v-u)^{2H-2}\text{ .}
	\]
	By Cauchy-Schwarz the same estimate as above holds in the case $u^+_\ell = v^-_\ell$, with a constant in the second inequality given by the fact that $v-u \leq 2\varrho_\ell$. Let $u^+_\ell < v^-_\ell$: we have, by \eqref{eq:requiredIneq} and for any $H'$ as in the statement
	\begin{align*}
			|\mathbb E[\dot X_u^\ell \dot X^\ell_v]| &= \varrho^{-2}_\ell \bigg|\int_{[u^-_\ell, u^+_\ell] \times [v^-_\ell, v^+_\ell]} \partial_{12}R(u,v) \dif u \dif v \bigg| \\
			&\lesssim \varrho^{-2}_\ell \int_{[u^-_\ell, u^+_\ell] \times [v^-_\ell, v^+_\ell]} (v-u)^{2H-2} \dif u \dif v  \\
			&\leq (v_\ell^- - u_\ell^+)^{2H-2} \\
			&\lesssim ((v^+_\ell-u^-_\ell) \wedge 1/2)^{2H'-2} \\
			&\leq ((v-u)\wedge 1/2)^{2H'-2} 
	\end{align*}
In the second-last inequality we have used that there exists some $L$ s.t.\ for all $\ell \geq L$
\[
\vartheta^{2H-2} \leq (\vartheta + 2\varrho_\ell)^{2H'-2}
\]
for all $\vartheta \in [\varrho_\ell, 1/2]$.

We now consider terms corresponding to maximal sequences of consecutive pairings, i.e.\
\begin{equation}\label{eq:n2k}
\int_{s < u_1 < v_1 < \ldots < u_k < v_k < t} \mathbb E[\dot X^\ell_{u_1} \dot X^\ell_{v_1}]\cdots \mathbb E[\dot X^\ell_{u_k} \dot X^\ell_{v_k}] \dif u_1 \dif v_1 \cdots \dif u_k \dif v_k\text{ .}
\end{equation}
It is always possible (e.g.\ by Kolmogorov's extension theorem) to add independent components to $X$. With this in mind, by Wick's theorem we may write the above integral as $\bbE \EuS(X^\ell)_{st}^{\alpha_1 \alpha_1 \ldots \alpha_k \alpha_k}$ with $\alpha_i \neq \alpha_j$ for all $i \neq j$. By the shuffle identity \eqref{eq:shuffle} we have
\begin{equation}\label{eq:consecInduc}
	\begin{split}
		&\sum_{h = 0}^{n} \EuS(X^\ell)_{st}^{\alpha_1\alpha_1 \ldots \alpha_h \alpha_h \beta\beta \alpha_{h+1} \alpha_{h+1} \ldots \alpha_k\alpha_k} \\
		={}&\EuS(X^\ell)_{st}^{\alpha_1 \alpha_1 \ldots \alpha_k \alpha_k} \EuS(X^\ell)_{st}^{\beta\beta} -  \sum_{0 \leq i < j \leq k} \EuS(X^\ell)_{st}^{\alpha_1 \alpha_1 \ldots \alpha_i \alpha_i \beta \alpha_{i+1} \alpha_{i+1} \ldots \alpha_j \alpha_j \beta \alpha_{j+1} \alpha_{j+1} \ldots \alpha_k \alpha_k} \\
		&-\sum_{0 \leq i < j \leq k} \EuS(X^\ell)_{st}^{\alpha_1 \alpha_1 \ldots \alpha_i \alpha_i \beta \alpha_{i+1} \alpha_{i+1} \ldots \alpha_{j+1} \beta \alpha_{j+1} \ldots \alpha_k \alpha_k} \\
		&-\sum_{0 \leq i < j \leq k} \EuS(X^\ell)_{st}^{\alpha_1 \alpha_1 \ldots \alpha_i \beta \alpha_i \ldots \alpha_j \alpha_j \beta \alpha_{j+1} \alpha_{j+1} \ldots \alpha_k \alpha_k} \\
		&- \sum_{0 \leq i < j \leq k}\EuS(X^\ell)_{st}^{\alpha_1 \alpha_1 \ldots \alpha_i \beta\alpha_i  \ldots \alpha_j \beta \alpha_j \ldots \alpha_k \alpha_k} \\
		&- \sum_{h = 0}^k  	\Big(\EuS(X^\ell)_{st}^{\alpha_1 \alpha_1 \ldots \alpha_h \beta\alpha_h \beta \alpha_{h+1}\alpha_{h+1} \ldots \alpha_k \alpha_k} + \EuS(X^\ell)_{st}^{\alpha_1 \alpha_1 \ldots \alpha_{h-1}\alpha_{h-1} \beta\alpha_h\beta \alpha_h \ldots \alpha_k \alpha_k} \Big)\text{ .}
	\end{split}
\end{equation}
When shuffling we have separated the cases in which all $\alpha_h \alpha_h$ and $\beta\beta$ occur as consecutive pairs, from those in which at least one such pair is separated. We now take expectations: note that both independence and equal distribution of components are used.
\begin{equation}\label{eq:consecInduc1}
	(k+1)\bbE \EuS(X^\ell)_{st}^{\alpha_1 \alpha_1 \ldots \alpha_k \alpha_k \beta \beta}	= \bbE \EuS(X^\ell)_{st}^{\alpha_1 \alpha_1 \ldots \alpha_k \alpha_k} \cdot \tfrac 12\bbE [(X^\ell)^2_{st}] - \sum_Q Q^\ell_{st}
\end{equation}
where we are summing over a finite number of diagrams $Q$ whose longest sequence of consecutive pairings contains $k$ pairs or fewer.

We now prove the statement in the case in which $P$ has no single nodes, by induction on $n$. For $n = 0$, $P^\ell_{st} = \varnothing^\ell_{st} \equiv 1$ there is nothing to show. Let $n \geq 1$, and assume we have rewritten the integral according to \eqref{eq:3ints} (where the middle integral may be skipped, since there are no Malliavin derivatives). If $P$ is not given by a sequence of $n/2$ consecutive pairs, all maximal sequences of consecutive pairs in $P$ consist of fewer than $n$ pairs, and that thus the inductive hypothesis applies to them: this means that for each such sequence $Q$ with $k$ pairs, $|Q^\ell_{uv}| \lesssim (v-u)^{2kH'}$. Using the bounds for the first two types of integrand derived in the first part of this proof, the statement for $P$ then follows from \autoref{prop:finite} applied in the modified case of \autoref{rem:finite} and with exponent $H'$. Assume now $n = 2(k+1)$ and let $P$ be given by the diagram consisting of $k+1$ consecutive pairs: the only thing needed to conclude the induction is the bound. This follows from \eqref{eq:consecInduc1} thanks to the inductive hypothesis and the boundedness statement of \autoref{lem:onePair}.

Finally, we consider the general case in which $P$ may have single nodes. This follows again by writing $P^\ell_{st}$ in nested form, bounding terms corresponding to non-consecutive pairings as done above, and bounding the middle integral in \eqref{eq:3ints} thanks to the boundedness statement of \autoref{lem:middle}. When invoking this lemma, $f_\ell$ is going to be a product of terms of the form \eqref{eq:n2k} (with the extrema $s$ and $t$ replaced with variables $u_i$ and $u_j$ already integrated in the outer or middle integral), which as already proved is bounded by $\lesssim (t-s)^{2H'k}$: this yields the required bound overall.
\end{proof}
\end{lem}
	
	\begin{prop}[Convergence]\label{prop:conv}
		The functions $[0,T]^m \to \bbR$ of \autoref{defn:seql} individually converge a.e.\ to those of \autoref{def:defIntegral}: for $P \in \mathcal P^n_m$ it holds that
		\begin{equation}\label{eq:convEq}
			P^\ell_{st} \xrightarrow{\ell \to \infty} P_{st}\text{ .}
		\end{equation}
	Moreover $|P_{st}| \lesssim |P|_{st}$ (the integrals of \autoref{prop:finite}) uniformly on $\Delta^m[s,t]$.
		\begin{proof}
			The inequality is an absolute estimate of $P_{st}$ using \eqref{eq:requiredIneq} and \eqref{eq:12REst}. The structure of the proof of the first statement closely mirrors that of the previous lemma: we first consider the case in which $P$ does not have single nodes. For $u^-_\ell < v^-_\ell$
			\begin{align*}
				\mathbb E[\dot X_u^\ell \dot X^\ell_v] = \varrho^{-2}_\ell R(\Delta(u^-_\ell,u^+_\ell),\Delta(v^-_\ell,v^+_\ell)) = \partial_{12}R(\overline u, \overline v)
			\end{align*}
			for some $\overline u \in (u^-_\ell,u^+_\ell)$, $\overline v \in (v^-_\ell,v^+_\ell)$,  by the intermediate value theorem applied twice. Pointwise convergence $\mathbb E[\dot X_u^\ell \dot X^\ell_v] \to \partial_{12}R(u,v)$ then holds by continuity of $\partial_{12}R$ and thanks to the fact that for any $u < v$ there exists $L$ s.t.\ $u^-_\ell < v^-_\ell$ for all $\ell \geq L$. This takes care of convergence of terms corresponding to non-consecutive pairings (of course, the same holds for consecutive pairings, but is not useful since $\partial_{12}R(u,v)$ may not be integrable in this case).
			
			We now proceed by induction on $n$. For $n = 0$ there is nothing to prove, so let $n \geq 1$ and first consider the case in which $P$ is not given by a sequence of $n/2$ consecutive pairs: the statement follows by dominated convergence applied to the outer integral in \eqref{eq:3ints}, by the above and the inductive hypothesis applied to sequences of consecutive nodes of length less than $n$, in conjunction with \autoref{lem:dom}. Let now $n = 2(k+1)$ and let $P$ be given by the diagram consisting of $k+1$ consecutive pairs: recalling the argument (and indexing notation) of the previous proof, we have $P_{st}^\ell = \bbE \EuS(X^\ell)_{st}^{\alpha_1 \alpha_1 \ldots \alpha_k \alpha_k \beta \beta}$, which is convergent since $\EuS(X^\ell)_{st} \to \EuS(X)_{st}$ in $L^2$. By the same calculation of \eqref{eq:consecInduc} applied to $X$ instead of to $X^\ell$, and taking expectations
			\begin{equation}\label{eq:limitPl}
				\begin{split}
				&(k+1)\lim_{\ell \to \infty} P_{st}^\ell\\
				={}&(k+1)\bbE \EuS(X)_{st}^{\alpha_1 \alpha_1 \ldots \alpha_k \alpha_k \beta \beta}  \\
				={} & \bbE \EuS(X)_{st}^{\alpha_1 \alpha_1 \ldots \alpha_k \alpha_k} \cdot \tfrac 12\bbE [(X)^2_{st}] - \sum_{0 \leq i < j \leq k} \bbE\EuS(X)_{st}^{\alpha_1 \alpha_1 \ldots \alpha_i \alpha_i \beta \alpha_{i+1} \alpha_{i+1} \ldots \alpha_j \alpha_j \beta \alpha_{j+1} \alpha_{j+1} \ldots \alpha_k \alpha_k} \\
				&-\sum_{0 \leq i < j \leq k} \bbE\EuS(X)_{st}^{\alpha_1 \alpha_1 \ldots \alpha_i \alpha_i \beta \alpha_{i+1} \alpha_{i+1} \ldots \alpha_{j+1} \beta \alpha_{j+1} \ldots \alpha_k \alpha_k} \\
				&-\sum_{0 \leq i < j \leq k} \bbE\EuS(X)_{st}^{\alpha_1 \alpha_1 \ldots \alpha_i \beta \alpha_i \ldots \alpha_j \alpha_j \beta \alpha_{j+1} \alpha_{j+1} \ldots \alpha_k \alpha_k} \\
				&- \sum_{0 \leq i < j \leq k}\bbE\EuS(X)_{st}^{\alpha_1 \alpha_1 \ldots \alpha_i \beta\alpha_i  \ldots \alpha_j \beta \alpha_j \ldots \alpha_k \alpha_k} \\
				&- \sum_{h = 0}^k \Big(\bbE\EuS(X)_{st}^{\alpha_1 \alpha_1 \ldots \alpha_h \beta\alpha_h \beta \alpha_{h+1}\alpha_{h+1} \ldots \alpha_k \alpha_k} + \bbE\EuS(X)_{st}^{\alpha_1 \alpha_1 \ldots \alpha_{h-1}\alpha_{h-1} \beta\alpha_h\beta \alpha_h \ldots \alpha_k \alpha_k} \Big)\text{ .}
				\end{split}
			\end{equation} 
		We now expand the product: by the inductive hypothesis and \autoref{lem:onePair}, and using Fubini's theorem we have (setting $u_0 =s = w_0$)
			\begin{align*}
			 &\bbE \EuS(X)_{st}^{\alpha_1 \alpha_1 \ldots \alpha_k \alpha_k} \cdot\bbE [X^2_{st}] \\
			={} & \int_{s < u_1 < \ldots < u_k < t} \big[ \tfrac 12 R'(u_1) - \partial_2R(u_0,u_1) \big] \cdots \big[ \tfrac 12 R'(u_k) - \partial_2R(u_{k-1},u_k) \big] \dif u_1 \cdots \dif u_k \\
			&\quad\cdot \int_s^t \big[\tfrac 12 R'(v) - \partial_2R(s,v) \big] \dif v \\
			={} & \int_{\substack{s < u_1 < \ldots < u_k <t \\ s < v < t}} \big[ \tfrac 12 R'(u_1) - \partial_2R(u_0,u_1) \big] \cdots \big[ \tfrac 12 R'(u_k) - \partial_2R(u_{k-1},u_k) \big] 	\\
			&\quad\cdot\big[\tfrac 12 R'(v) - \partial_2R(s,v) \big] \dif u_1 \cdots \dif u_k \dif v \\
	={}&\sum_{j = 0}^k\int_{s < u_1 < \ldots < u_j < v < u_{j+1} < \ldots < u_k < t}\big[ \tfrac 12 R'(u_1) - \partial_2R(u_0,u_1) \big] \cdots \big[ \tfrac 12 R'(u_k) - \partial_2R(u_{k-1},u_k) \big] \\
	&\qquad\cdot\big[\tfrac 12 R'(v) - \partial_2R(s,v) \big] \dif u_1 \cdots \dif u_k \dif v\text{ .}
			\end{align*}
	Note that the use of Fubini's theorem is justified in view of \eqref{eq:12REst} applied to absolutely bound each integral above, and \autoref{prop:finite}. Writing 
	\[
	\partial_2R(\Delta(x,y),z) \coloneqq \partial_2R(y,z) - \partial_2R(x,z) = \int_x^y \partial_{12}R(w,z) \dif w \text{ ,}
	\] we expand each summand:
	\begin{align*}
		&\int_{s < u_1 < \ldots < u_j < v < u_{j+1} < \ldots < u_k < t}\big[ \tfrac 12 R'(u_1) - \partial_2R(u_0,u_1) \big] \cdots \big[ \tfrac 12 R'(u_k) - \partial_2R(u_{k-1},u_k) \big] \\
		&\quad\cdot\big[\tfrac 12 R'(v) - \partial_2R(s,v) \big] \dif u_1 \cdots \dif u_k \dif v \\
		={}&\int_{s < u_1 < \ldots < u_j < v < u_{j+1} < \ldots < u_k < t}\big[ \tfrac 12 R'(u_1) - \partial_2R(u_0,u_1) \big] \cdots \big[ \tfrac 12 R'(u_j) - \partial_2R(u_{j-1},u_j) \big] \\
		&\quad\cdot \big[\tfrac 12 R'(v) - \partial_2 R(u_j,v) + \textstyle \sum_{i=0}^{j-1}\partial_2R(\Delta(u_i,u_{i+1}),v) \big]\\
		&\quad\cdot\big[ \tfrac 12 R'(u_{j+1}) - \partial_2R(v,u_{j+1}) + \partial_2R(\Delta(u_j,v),u_{j+1}) \big]\\
		&\quad \cdot \big[\tfrac 12 R'(u_{j+2}) - \partial_2R(u_{j+1},u_{j+2}) \big] \cdots \big[\tfrac 12 R'(u_k) - \partial_2R(u_{k-1},u_k) \big]  \dif u_1 \cdots \dif u_k \dif v \\
		={}&\int_{s < u_1 < \ldots < u_j < v < u_{j+1} < \ldots < u_k < t} \big[ \tfrac 12 R'(u_1) - \partial_2R(u_0,u_1) \big] \cdots \big[\tfrac 12 R'(v) - \partial_2 R(u_j,v)\big] \\
		&\quad\cdot\big[\tfrac 12 R'(u_{j+1}) - \partial_2R(u_{j+1},v) \big] \cdots \big[\tfrac 12 R'(u_k) - \partial_2R(u_{k-1},u_k) \big]  \dif u_1 \cdots \dif u_k \dif v \\
		&+\int_{s < u_1 < \ldots < u_j < r < v < u_{j+1} < \ldots < u_k < t} \big[ \tfrac 12 R'(u_1) - \partial_2R(u_0,u_1) \big] \cdots \big[ \tfrac 12 R'(u_j) - \partial_2R(u_{j-1},u_j) \big] \\
		&\quad\cdot \big[\tfrac 12 R'(v) - \partial_2 R(r,v) + \partial_2R(\Delta(u_j,r),v)\big]\partial_{12}R(r,u_{j+1})  \dif u_1 \cdots \dif u_k \dif r \dif v \\
		&+\sum_{i = 0}^{j-1} \int_{s < u_1 < \ldots < u_i < q < u_{i+1} < \ldots < u_j < v < u_{j+1} < \ldots < u_k < t}  \big[ \tfrac 12 R'(u_1) - \partial_2R(u_0,u_1) \big] \\
		&\qquad \cdots  \partial_{12}R(q,v)\big[ \tfrac 12 R'(u_{i+1}) -\partial_2R(q,u_{i+1}) + \partial_2R(\Delta(u_i,q),u_{i+1}) \big] \cdots \big[\tfrac 12 R'(u_{j+1}) - \partial_2R(v,u_{j+1}) \big] \\
		&\qquad\cdots \big[ \tfrac 12 R'(u_k) - \partial_2R(u_{k-1},u_k) \big] \dif u_1 \cdots \dif u_k \dif q \dif v \\
		&+\sum_{i = 0}^{j-1} \int_{s < u_1 < \ldots < u_i < q < u_{i+1} < \ldots < u_j < r < v < u_{j+1} < \ldots < u_k < t}  \big[ \tfrac 12 R'(u_1) - \partial_2R(u_0,u_1) \big] \\
		&\qquad\cdots  \partial_{12}R(q,v)\big[ \tfrac 12 R'(u_{i+1}) -\partial_2R(q,u_{i+1}) + \partial_2R(\Delta(u_i,q),u_{i+1}) \big] \cdots \partial_{12}R(r,u_{j+1}) \\
		&\qquad\cdots \big[  \tfrac 12 R'(u_k) - \partial_2R(u_{k-1},u_k) \big]\dif u_1 \cdots \dif u_k \dif q \dif r \dif v\\
		={}& (\,\tikz[baseline = -0.1ex]{
			\draw [decorate,decoration={brace,amplitude=4pt,mirror,raise=0.6ex}]
			(-0.1 ,0) -- (1,0);
			\draw[fill] (0,0) circle [radius=0.04];
			\draw[fill] (0.9,0) circle [radius=0.04];
			\draw (0,0)-- (0,0.275);
			\draw (0.9,0)-- (0.9,0.275);
			\node at (0.45,0){$\ldots$};
			\node at (0.45,-0.4){$\scriptstyle k+1$};
		}\,)_{st} + (\,\tikz[baseline = -0.1ex]{
		\draw [decorate,decoration={brace,amplitude=4pt,mirror,raise=0.6ex}]
		(-0.1 ,0) -- (1,0);
		\draw [decorate,decoration={brace,amplitude=4pt,mirror,raise=0.6ex}]
		(2,0) -- (3.1,0);
		\draw[fill] (0,0) circle [radius=0.04];
		\draw[fill] (0.9,0) circle [radius=0.04];
		\draw (0,0)-- (0,0.275);
		\draw (0.9,0)-- (0.9,0.275);
		\node at (0.45,0){$\ldots$};
		\node at (0.45,-0.4){$\scriptstyle j$};
		\draw[fill] (1.2,0) circle [radius=0.04];
		\draw[fill] (1.5,0) circle [radius=0.04];
		\draw[fill] (1.8,0) circle [radius=0.04];
		\draw[fill] (2.1,0) circle [radius=0.04];
		\draw[fill] (3,0) circle [radius=0.04];
		\draw (2.1,0)-- (2.1,0.275);
		\draw (3,0)-- (3,0.275);
		\draw (1.5,0)-- (1.5,0.275);
		\node at (2.55,0){$\ldots$};
		\node at (2.55,-0.4){$\scriptstyle k-j-1$};
		\draw (1.2,0) .. controls (1.2,0.5) and (1.8,0.5) .. (1.8,0);
	}\,)_{st} + (\,\tikz[baseline = -0.1ex]{
	\draw [decorate,decoration={brace,amplitude=4pt,mirror,raise=0.6ex}]
	(-0.1 ,0) -- (1,0);
	\draw [decorate,decoration={brace,amplitude=4pt,mirror,raise=0.6ex}]
	(2.3,0) -- (3.4,0);
	\draw[fill] (0,0) circle [radius=0.04];
	\draw[fill] (0.9,0) circle [radius=0.04];
	\draw (0,0)-- (0,0.275);
	\draw (0.9,0)-- (0.9,0.275);
	\node at (0.45,0){$\ldots$};
	\node at (0.45,-0.4){$\scriptstyle j$};
	\draw[fill] (1.2,0) circle [radius=0.04];
	\draw[fill] (1.5,0) circle [radius=0.04];
	\draw[fill] (1.8,0) circle [radius=0.04];
	\draw[fill] (2.4,0) circle [radius=0.04];
	\draw[fill] (3.3,0) circle [radius=0.04];
	\draw[fill] (2.1,0) circle [radius=0.04];
	\draw (2.4,0)-- (2.4,0.275);
	\draw (3.3,0)-- (3.3,0.275);
	\node at (2.85,0){$\ldots$};
	\node at (2.85,-0.4){$\scriptstyle k-j-1$};
	\draw (1.2,0) .. controls (1.2,0.5) and (1.8,0.5) .. (1.8,0);
		\draw (1.5,0) .. controls (1.5,0.5) and (2.1,0.5) .. (2.1,0);
}\,)_{st} \\ 
&+ \sum_{i = 0}^{j-1} (\,\tikz[baseline = -0.1ex]{
	\draw [decorate,decoration={brace,amplitude=4pt,mirror,raise=0.6ex}]
	(-0.1 ,0) -- (1,0);
	\draw [decorate,decoration={brace,amplitude=4pt,mirror,raise=0.6ex}]
	(1.4 ,0) -- (2.5,0);
	\draw[fill] (0,0) circle [radius=0.04];
	\draw[fill] (0.9,0) circle [radius=0.04];
	\draw (0,0)-- (0,0.275);
	\draw (0.9,0)-- (0.9,0.275);
	\node at (0.45,0){$\ldots$};
	\node at (0.45,-0.4){$\scriptstyle i$};
	\draw[fill] (1.2,0) circle [radius=0.04];
	\draw[fill] (1.5,0) circle [radius=0.04];
	\draw[fill] (2.4,0) circle [radius=0.04];
	\draw (1.5,0)-- (1.5,0.275);
	\draw (2.4,0)-- (2.4,0.275);
	\draw[fill] (2.7,0) circle [radius=0.04];
	\node at (1.95,0){$\ldots$};
	\node at (1.95,-0.4){$\scriptstyle j-i$};
	\draw (1.2,0) .. controls (1.1,0.6) and (2.8,0.6) .. (2.7,0);
	\draw [decorate,decoration={brace,amplitude=4pt,mirror,raise=0.6ex}]
	(2.9,0) -- (4,0);
	\draw[fill] (3,0) circle [radius=0.04];
	\draw[fill] (3.9,0) circle [radius=0.04];
	\draw (3,0)-- (3,0.275);
	\draw (3.9,0)-- (3.9,0.275);
	\node at (3.45,0){$\ldots$};
	\node at (3.45,-0.4){$\scriptstyle k-j$};
}\,)_{st}  + \sum_{i = 0}^{j-1} (\,\tikz[baseline = -0.1ex]{
\draw [decorate,decoration={brace,amplitude=4pt,mirror,raise=0.6ex}]
(-0.1 ,0) -- (1,0);
\draw[fill] (0,0) circle [radius=0.04];
\draw[fill] (0.9,0) circle [radius=0.04];
\draw (0,0)-- (0,0.275);
\draw (0.9,0)-- (0.9,0.275);
\node at (0.45,0){$\ldots$};
\node at (0.45,-0.4){$\scriptstyle i$};
\draw[fill] (1.2,0) circle [radius=0.04];
\draw[fill] (1.5,0) circle [radius=0.04];
\draw[fill] (1.8,0) circle [radius=0.04];
\draw (1.2,0) .. controls (1.2,0.5) and (1.8,0.5) .. (1.8,0);
\draw [decorate,decoration={brace,amplitude=4pt,mirror,raise=0.6ex}]
(2,0) -- (3.1,0);
\draw[fill] (2.1,0) circle [radius=0.04];
\draw[fill] (3,0) circle [radius=0.04];
\draw (2.1,0)-- (2.1,0.275);
\draw (3,0)-- (3,0.275);
\node at (2.55,0){$\ldots$};
\node at (2.55,-0.4){$\scriptstyle j-i-1$};
\draw[fill] (3.3,0) circle [radius=0.04];
\draw (1.2+0.3,0) .. controls (1.1+0.3,0.6) and (3.1+0.3,0.6) .. (3+0.3,0);
\draw [decorate,decoration={brace,amplitude=4pt,mirror,raise=0.6ex}]
(-0.1+3.6 ,0) -- (1+3.6,0);
\draw[fill] (0+3.6,0) circle [radius=0.04];
\draw[fill] (0.9+3.6,0) circle [radius=0.04];
\draw (0+3.6,0)-- (0+3.6,0.275);
\draw (0.9+3.6,0)-- (0.9+3.6,0.275);
\node at (0.45+3.6,0){$\ldots$};
\node at (0.45+3.6,-0.4){$\scriptstyle k-j$};
} \,)_{st}\\
&+\sum_{i = 0}^{j-1} (\,\tikz[baseline = -0.1ex]{
\draw [decorate,decoration={brace,amplitude=4pt,mirror,raise=0.6ex}]
(-0.1 ,0) -- (1,0);
	\draw [decorate,decoration={brace,amplitude=4pt,mirror,raise=0.6ex}]
(1.4 ,0) -- (2.5,0);
\draw[fill] (0,0) circle [radius=0.04];
\draw[fill] (0.9,0) circle [radius=0.04];
\draw (0,0)-- (0,0.275);
\draw (0.9,0)-- (0.9,0.275);
\node at (0.45,0){$\ldots$};
\node at (0.45,-0.4){$\scriptstyle i$};
\draw[fill] (1.2,0) circle [radius=0.04];
\draw[fill] (1.5,0) circle [radius=0.04];
\draw[fill] (2.4,0) circle [radius=0.04];
\draw (1.5,0)-- (1.5,0.275);
\draw (2.4,0)-- (2.4,0.275);
\node at (1.95,0){$\ldots$};
\draw[fill] (2.7,0) circle [radius=0.04];
\draw[fill] (3,0) circle [radius=0.04];
\draw[fill] (3.3,0) circle [radius=0.04];
	\node at (1.95,-0.4){$\scriptstyle j-i$};
\draw (1.2,0) .. controls (1.1,0.6) and (3.1,0.6) .. (3,0);
\draw (2.7,0) .. controls (2.7,0.5) and (3.3,0.5) .. (3.3,0);
\draw [decorate,decoration={brace,amplitude=4pt,mirror,raise=0.6ex}]
(2.9+0.6,0) -- (4+0.6,0);
\draw[fill] (3+0.6,0) circle [radius=0.04];
\draw[fill] (3.9+0.6,0) circle [radius=0.04];
\draw (3+0.6,0)-- (3+0.6,0.275);
\draw (3.9+0.6,0)-- (3.9+0.6,0.275);
\node at (3.45+0.6,0){$\ldots$};
\node at (3.45+0.6,-0.4){$\scriptstyle k-j-1$};
}\,)_{st} + \sum_{i = 0}^{j-1} (\,\tikz[baseline = -0.1ex]{
\draw [decorate,decoration={brace,amplitude=4pt,mirror,raise=0.6ex}]
(-0.1 ,0) -- (1,0);
\draw[fill] (0,0) circle [radius=0.04];
\draw[fill] (0.9,0) circle [radius=0.04];
\draw (0,0)-- (0,0.275);
\draw (0.9,0)-- (0.9,0.275);
\node at (0.45,0){$\ldots$};
\node at (0.45,-0.4){$\scriptstyle i$};
\draw[fill] (1.2,0) circle [radius=0.04];
\draw[fill] (1.5,0) circle [radius=0.04];
\draw[fill] (1.8,0) circle [radius=0.04];
\draw (1.2,0) .. controls (1.2,0.5) and (1.8,0.5) .. (1.8,0);
\draw (1.2+2.1,0) .. controls (1.2+2.1,0.5) and (1.8+2.1,0.5) .. (1.8+2.1,0);
\draw [decorate,decoration={brace,amplitude=4pt,mirror,raise=0.6ex}]
(2,0) -- (3.1,0);
\draw[fill] (2.1,0) circle [radius=0.04];
\draw[fill] (3,0) circle [radius=0.04];
\draw (2.1,0)-- (2.1,0.275);
\draw (3,0)-- (3,0.275);
\node at (2.55,0){$\ldots$};
\node at (4.65,0){$\ldots$};
\node at (2.55,-0.4){$\scriptstyle j-i-1$};
\node at (4.65,-0.4){$\scriptstyle k-j-2$};
\draw[fill] (3.3,0) circle [radius=0.04];
\draw (1.2+0.3,0) .. controls (1.1+0.3,0.6) and (3.1+0.6,0.6) .. (3+0.6,0);
\draw [decorate,decoration={brace,amplitude=4pt,mirror,raise=0.6ex}]
(-0.1+3.6+0.6 ,0) -- (1+3.6+0.6,0);
\draw[fill] (0+3.6+0.6,0) circle [radius=0.04];
\draw[fill] (0.9+3.3,0) circle [radius=0.04];
\draw[fill] (0.9+3,0) circle [radius=0.04];
\draw[fill] (5.1,0) circle [radius=0.04];
\draw[fill] (3.6,0) circle [radius=0.04];
\draw (0+3.6+0.6,0)-- (0+3.6+0.6,0.275);
\draw (0.9+3.6+0.6,0)-- (0.9+3.6+0.6,0.275);
} \,)_{st} \\
={}& (k+1)P_{st} + \bbE \EuS(X)_{st}^{\alpha_1 \alpha_1 \ldots \alpha_j\alpha_j \beta \alpha_{j+1}\alpha_{j+1} \beta \alpha_{j+2} \alpha_{j+2} \ldots \alpha_k \alpha_k} + \bbE \EuS(X)^{\alpha_1 \alpha_1 \ldots \alpha_j\alpha_j \beta \alpha_{j+1}\beta\alpha_{j+1} \ldots \alpha_k \alpha_k}_{st} \\
&+ \sum_{i = 0}^{j-1}  \EuS(X)_{st}^{\alpha_1 \alpha_1 \ldots \alpha_i\alpha_i \beta \alpha_{i+1}\alpha_{i+1} \ldots \alpha_j \alpha_j \beta \alpha_{j+1} \alpha_{j+1} \ldots \alpha_k \alpha_k} + \sum_{i = 0}^{j-1}  \EuS(X)_{st}^{\alpha_1 \alpha_1 \ldots \alpha_{i+1}\beta \alpha_{i+1} \ldots \alpha_j \alpha_j \beta \alpha_{j+1} \alpha_{j+1} \ldots \alpha_k \alpha_k} \\
&+\sum_{i = 0}^{j-1}  \EuS(X)_{st}^{\alpha_1 \alpha_1 \ldots \alpha_i \alpha_i\beta \alpha_{i+1}\alpha_{i+1} \ldots \alpha_{j+1} \beta \alpha_{j+1} \ldots \alpha_k \alpha_k} +\sum_{i = 0}^{j-1}  \EuS(X)_{st}^{\alpha_1 \alpha_1 \ldots \alpha_{i+1} \beta \alpha_{i+1} \ldots \alpha_{j+1} \beta \alpha_{j+1} \ldots \alpha_k \alpha_k} \text{ .}
\end{align*}
It now follows by substitution into the sum $\sum_{j = 0}^k$ and simplifying in \eqref{eq:limitPl} that $\lim_\ell P^\ell_{st} = P_{st}$.

Finally, we consider diagrams that contain single nodes. In order to invoke \autoref{lem:middle} we must argue that $f_\ell \to f$ uniformly (uniform boundedness holds by the previous lemma). This again follows from the fact that $f_\ell$ can be written as a product of expected signatures of $X^\ell$, each of which converges uniformly in $\ell$ as a function of its extrema: recalling the notations for truncation and projection introduced in \autoref{sec:back} and the definition of inhomogeneous $p$-variation distance \cite[\S 8.1.2]{FV10}, we have
\begin{align*}
\sup_{u < v}|\mathbb E \EuS(X)_{uv}^{(n)} - \mathbb E \EuS(X^\ell)_{uv}^{(n)}| &\leq \mathbb E \sup_{u < v} |\EuS(X)_{uv}^{(n)}-  \EuS(X^\ell)_{uv}^{(n)}| \\
&\leq \lVert \rho_{p\mhyphen\mathrm{var}} (\EuS^\p(X^\ell), \EuS^\p(X)) \rVert_{L^1} \\
&\xrightarrow{\ell \to \infty} 0
\end{align*}
for $p > (1/H) \vee n$, where we have used \cite[Theorem 1]{FR14}. The statement now follows once again by dominated convergence and Fubini's theorem.
\end{proof}			
\end{prop}

We are ready to put it all together:
\begin{proof}[Proof of \autoref{thm:main}]
	\begin{align}
		\mathscr w^m \EuS(X)_{st}^{\gamma_1,\ldots,\gamma_n} &= \frac{1}{m!}\skoo^m \big( \bbE \EuD^m \EuS(X)_{st}^{\gamma_1,\ldots,\gamma_n}\label{eq:proof1} \big) \\
		&= \frac{1}{m!}\skoo^m \lim_{\ell \to \infty} \big( \bbE \EuD^m \EuS(X^\ell)_{st}^{\gamma_1,\ldots,\gamma_n} \big) \label{eq:proof2}\\
		&= \skoo^m \sum_{P \in \mathcal P^n_m} \lim_{\ell \to \infty}P_{st}^{\ell;\gamma_1,\ldots,\gamma_n} \label{eq:proof3}\\
		&= \sum_{P \in \mathcal P^n_m}  \skoo^m P^{\gamma_1,\ldots,\gamma_n}_{st}\label{eq:proof4}
	\end{align}
	In \eqref{eq:proof1} we have used Stroock's formula \eqref{eq:stroock}, which is possible since $\EuS(X)_{st}^{\gamma_1,\ldots,\gamma_n} \in \mathbb D^{\infty,2}$: this is because $\EuS(X^\ell)_{st} \to \EuS(X)_{st}$ in $\bigoplus_{k \leq n}\mathscr W^k$ which is closed in $L^2 \Omega$. In \eqref{eq:proof2} we have used that convergence of $\EuS(X^\ell)_{st}$ actually holds in $\mathbb D^{\infty,2}$, since the norm of $\mathbb D^{\infty,2}$ is dominated by the $L^2$ norm in bounded Wiener chaos \cite[Proposition 1.2.2]{Nua06}. \eqref{eq:proof3} uses \autoref{lem:expMal} and \eqref{eq:proof4} is just the statement (required by our definition of membership of a function to $\mathcal H^{\otimes m}$ \autoref{def:fH}) that $P_{st}^{\ell;\gamma_1,\ldots,\gamma_n}$ converges a.e.\ boundedly to $P^\ell_{st}$, which holds by \autoref{prop:conv} and \autoref{lem:dom}. As argued in the previous two proofs, $P_{st}^{\gamma_1,\ldots,\gamma_n}$ can always be expressed as the expected signature evaluated on a word, up to augmenting $X$ with independent copies of itself: this can be used to infer that each $P_{st}^{\gamma_1,\ldots,\gamma_n}$ --- not just their sum --- belongs to $\mathbb D^{\infty,2}(\mathcal H^{\otimes m})$. This concludes the proof of the main result.
\end{proof}

\starsection{Conclusions and further directions}\addcontentsline{toc}{section}{Conclusions and further directions}\label{sec:concl}

By providing a single formula for the expected signature of fractional Brownian motion that holds for any Hurst parameter $H \in (1/4,1)$, this article closes a gap in the literature left open by \cite{Bau07}. Along the way, we have had opportunity to consider numerous other aspects of our computation, such as similar formulae for higher levels of the Wiener chaos expansion of the signature, and other examples of Gaussian processes.

We believe this work recommends a variety of applications and further investigations. First and foremost, it would be interesting to write stochastic Taylor expansions as suggested by \autoref{expl:taylor0}, under precise conditions on the vector fields, and by providing bounds on the mean square error. Making this calculation rigorous and providing precise asymptotic estimates such as those in \cite{Pa20} would be an interesting result, which could be applied to approximation problems for Gaussian RDEs on manifolds such as those considered in \cite{ABR19} for SDEs (although for this precise problem, the joint signature $\EuS(X,t)$ would have to be considered). A further step would involve proving conditional versions of the results in this paper, which would make it possible to estimate the error generated by multiple steps in an Euler scheme.

The fact that (e.g.\ for fBm) the integral $\bbE \EuS(X)^{\alpha_1\alpha_1 \cdots \alpha_k\alpha_k}_{st}$ with $\alpha_i \neq \alpha_j$ is actually convergent for any $H > 0$ raises the question of whether something can be said about the sequence $\EuS(X^\ell)^{\alpha_1\alpha_1 \cdots \alpha_k\alpha_k}_{st}$, i.e.\ by considering the particular word on which $\EuS(X^\ell)$, which is not convergent in probability for $H \leq 1/4$, is evaluated.

It would be interesting to express the expected signature of a Gaussian process as the exponential of a formal series of tensors, thus computing its signature cumulants \cite{BO20}: this is how the expected signature of Brownian motion \eqref{eq:sumSquares} is usually presented (with the series a finite sum), but the analogous formulation for Gaussian processes that are not martingales appears more difficult to write down.

A more computational goal, though not one that appears trivial, is to explicitly compute \autoref{thm:main} for certain semimartingales, such as the Brownian bridge, for which the integrals are all analytically solvable. An interesting question is how the relationship between Brownian motion and Brownian bridge is reflected by their expected signatures. It would also be helpful to see whether similar formulae to ours can be made available for non-centred Gaussian processes, e.g.\ general Ornstein-Uhlenbeck processes. Finally, it would be interesting to try to apply the main theorem to the Riemann-Liouville process \autoref{expl:RV}.

\appendix

\section{Equivalence with \cite{Bau07} for the expected signature of fBm at level $4$}\label{appendix}

In \cite[Theorem 34]{Bau07} the authors check that the explicit integral expressions for $\mathbb E \EuS (X)^{(n)}_{01}$ with $n = 2,4$, previously derived for $X$ a fractional Brownian motion of Hurst parameter $H \in (1/2 , 1)$, continue to be valid for $H \in (1/4 , 1)$. This is done by performing transformations on $\mathbb E \EuS(X^\ell)_{01}$ before passing to the limit. This calculation is specific to levels $2$ and $4$ and to fBm, and for this reason the expression for the expected signature is not immediately comparable to that obtained as a special case of \autoref{thm:main}. We devote this appendix to checking that the two agree.

Level $2$ is easy to check, since \eqref{eq:ES2} reduces to $\text{\textdelta}^{\alpha\beta}/2$. Referring to \autoref{expl:EX4}, we consider level $4$ using \eqref{eq:t-s}: starting with the first integral above, we have
\begin{align*}
	&\int_{0 < u < v < 1} \big[ \tfrac 12 R(\dif u) - R(s,\dif u) \big] \big[ \tfrac 12 R(\dif v) - R(u,\dif v) \big] \dif u \dif v \\
	={}&H^2 \int_{0 < u < v < 1} u^{2H-1} (v-u)^{2H-1} \dif u \dif v \\
	={}&H^2 \int_0^1 u^{2H-1} \bigg[\frac{(v-u)^{2H-1}}{2H}\bigg]_{u = 0}^1 \dif u \\
	={}&\frac{H^2}{2} \int_0^1 u^{2H-1} \bigg[\frac{(v-u)^{2H}}{2H}\bigg]_{v = u}^1 \dif u \\
	={}&\frac{H}{2} \int_0^1 u^{2H-1} (1-u)^{2H} \dif u \\
	={}&\frac{H}{4} \int_0^1 u^{2H-1} (1-u)^{2H-1} \dif u
\end{align*}
where the last identity can be verified by showing that the difference of the two integrands is odd about the point $u = 1/2$, which in turn is seen by observing that
\[
\frac{H}{4} u^{2H-1} (1-u)^{2H-1} - \frac{H}{2} u^{2H-1} (1-u)^{2H} + \frac{H}{4} (1-u)^{2H-1} u^{2H-1} - \frac{H}{2} (1-u)^{2H-1} u^{2H}
\]
has zero derivative and vanishes at $u = 1/2$. This shows equality with \cite[coefficient of the first term of $\Gamma^2_H$ in Corollary 33]{Bau07}. We proceed with the second integral in \autoref{expl:EX4}:
\begin{align*}
	&\int_{s < u < v < w < t} R(\dif u,\dif w) \big[ \tfrac 12 R(\dif v) - R(u,\dif v) \big] \dif u \dif v \dif w \\
	={} &H^2(2H -1)\int_{0 < u < v < w < 1} (w-u)^{2H-2} (v-u)^{2H-1} \dif u \dif v \dif w \\
	={} &H^2\int_{0 < u < v < 1} \Big[(1-u)^{2H-1}(v-u)^{2H-1} - (v-u)^{4H-2} \Big]  \dif u \dif v \\
	={} &\bigg(\frac{H}{2} -  \frac{H^2}{4H-1} \bigg)\int_0^1 (1-u)^{4H-1}\dif u \\
	={} &\frac{2H - 1}{8(4H-1)}\text{ .}
\end{align*}
For the third integral we have
\begin{align*}
	&\int_{0 < u < v < w < z < 1} R(\dif u, \dif w) R(\dif v, \dif z)\\
	={} &H^2(2H - 1)^2 \int_{0 < u < v < w < z < 1} (w-u)^{2H-2}(z-v)^{2H-2} \dif u \dif v \dif w \dif z\\
	={} &H^2(2H - 1) \int_{0 < u < v < z < 1} \Big[ (z-u)^{2H-1} (z-v)^{2H-2}  - (v-u)^{2H-1} (z-v)^{2H-2} \Big] \dif u \dif v \dif z \\
	={} &\frac{H(2H - 1)}{2} \int_{0 < v < z < 1} \Big[  z^{2H} (z-v)^{2H-2} - (z-v)^{4H-2} - v^{2H} (z-v)^{2H-2}  \Big]  \dif v \dif z \\
	={} &\frac{H(2H - 1)}{2}\int_{0 < v < z < 1}(z^{2H}-v^{2H}) (z-v)^{2H-2} \dif v \dif z - \frac{H(2H - 1)}{4H-1}\int_0^1 (1-v)^{4H-1} \dif v  \\
	={} &\frac{H(2H - 1)}{2}\int_{0 < v < z < 1}(z^{2H}-v^{2H}) (z-v)^{2H-2} \dif v \dif z - \frac{2H - 1}{4(4H-1)} \\
	={} &\frac{H}{2}\int_0^1(1-v^{2H}) (1-v)^{2H-1} \dif v - H^2 \int_{0 < v < z < 1} z^{2H-1} (z-v)^{2H-1} - \frac{2H - 1}{4(4H-1)} \\
	={} &\frac{H}{4(4H-1)} - \frac{H}{4}\int_0^1 v^{2H-1} (1-v)^{2H-1} \dif v\text{ .}
\end{align*}
In the integration by parts we have used that $\lim_{z\to v^+} (z^{2H} - v^{2H})(z-v)^{2H-1} = 0$ which can be shown by using that for $1/4 < H < 1/2$
\[
0 \leq (z^{2H} - v^{2H})(z-v)^{2H-1} \leq (z-v)^{4H-1} \xrightarrow{z\to v^+} 0
\]
since $z^{2H}-v^{2H} < (z-v)^{2H}$ for $0 < v < z$ and $H < 1/2$. In the last identity we have used a similar symmetry argument as the one used in the first calculation, solved trivial integrals and rearranged terms. Note how this calculation would have been simpler if $H \geq 1/2$ since it would not have been necessary to integrate by parts to avoid integrating $(z-v)^{2H-2}$ (cf.\ \cite[Lemma 32]{Bau07}).

\subsection*{Data availability statement}
The research presented in this manuscript is theoretical in nature and does not refer to empirical data. All necessary information is contained in the cited references.

\subsection*{Author contribution statement}
Thomas Cass and Emilio Ferrucci are both authors of this paper. Emilio Ferrucci is the corresponding author.

\subsection*{Conflict of interest statement}
The authors declare no conflicts of interest.

\bibliographystyle{abbrv} 
\renewcommand\bibname{\sc References}
\bibliography{merged}\addcontentsline{toc}{chapter}{References}

\end{document}